\documentclass[11pt]{amsart}

\usepackage{amssymb,amsmath,amsthm}
\usepackage{color}
\usepackage{thmtools}
\usepackage{thm-restate}
\usepackage[shortlabels]{enumitem}
\usepackage{mathrsfs}
\usepackage[utf8]{inputenc}
\usepackage[T1]{fontenc}
\usepackage{dsfont}
\usepackage{soul}
\usepackage[margin=1in]{geometry}
\usepackage{mathtools}

\newtheorem{theorem}{Theorem}[section]
\newtheorem*{theorem*}{Theorem}
\newtheorem{lemma}[theorem]{Lemma}
\newtheorem{proposition}[theorem]{Proposition}

\newtheorem{question}[theorem]{Question}
\newtheorem{corollary}[theorem]{Corollary}
\newtheorem*{corollary*}{Corollary}

\theoremstyle{definition}
\newtheorem{definition}[theorem]{Definition}

\theoremstyle{remark}
\newtheorem{remark}[theorem]{Remark}

\numberwithin{equation}{section}

\newcommand{\frakc}{\mathfrak{c}}

\newcommand{\frakp}{\mathfrak{p}}

\newcommand{\fraks}{\mathfrak{s}}

\newcommand{\frakgr}{\mathfrak{gr}}

\newcommand{\eps}{\varepsilon}

\newcommand{\N}{\mathbb{N}}
\newcommand{\R}{\mathbb{R}}

\newcommand{\PP}{\mathbb{P}}

\newcommand{\aA}{\mathcal{A}}
\newcommand{\bB}{\mathcal{B}}
\newcommand{\cC}{\mathcal{C}}
\newcommand{\dD}{\mathcal{D}}
\newcommand{\eE}{\mathcal{E}}
\newcommand{\fF}{\mathcal{F}}

\newcommand{\mM}{\mathcal{M}}

\newcommand{\pP}{\mathcal{P}}

\DeclareMathOperator{\ba}{ba}

\DeclareMathOperator{\cf}{cf}

\DeclareMathOperator{\cov}{cov}

\newcommand{\ol}{\overline}

\newcommand{\rstr}{\restriction}

\newcommand{\sm}{\setminus}
\newcommand{\sub}{\subseteq}

\newcommand{\seq}[2]{\big\langle#1\colon\ #2\big\rangle}
\newcommand{\seqn}[1]{\big\langle#1\colon\ n\io\big\rangle}
\newcommand{\seqk}[1]{\big\langle#1\colon\ k\io\big\rangle}
\newcommand{\seql}[1]{\big\langle#1\colon\ l\io\big\rangle}
\newcommand{\seqm}[1]{\big\langle#1\colon\ m\io\big\rangle}
\newcommand{\seqp}[1]{\big\langle#1\colon\ p\io\big\rangle}
\newcommand{\seqq}[1]{\big\langle#1\colon\ q\io\big\rangle}

\newcommand{\iA}{\in\aA}

\newcommand{\io}{\in\mathbb{N}}

\newcommand{\elli}{{\ell_\infty}}

\newcommand{\Cantor}{{2^{\N}}}

\newcommand{\noproof}{\hfill$\Box$}

\begin{document}

\title[Algebra with the Grothendieck property but without the Nikodym property]{Construction under Martin's axiom of a Boolean algebra with the Grothendieck property but without the Nikodym property}
\author[D.\ Sobota]{Damian Sobota}
\address{Universit\"at Wien, Institut f\"ur Mathematik,
Kurt G\"odel Research Center, 
Wien, Austria}
\email{ein.damian.sobota@gmail.com}
\urladdr{www.logic.univie.ac.at/~{}dsobota}
\author[L.\ Zdomskyy]{Lyubomyr Zdomskyy}
\address{Technische Universität Wien, Institut für Diskrete Mathematik und Geometrie, Wien, Austria}
\email{lzdomsky@gmail.com}
\urladdr{dmg.tuwien.ac.at/zdomskyy}
\thanks{D. Sobota was supported by the Austrian Science Fund (FWF), Grant ESP 108-N. L. Zdomskyy was supported by the Austrian Science Fund (FWF), Grant I 5930-N}

\begin{abstract}
Improving a result of M. Talagrand, under the assumption of a weak form of Martin's axiom, we construct a totally disconnected compact Hausdorff space $K$ such that the Banach space $C(K)$ of continuous real-valued functions on $K$ is a Grothendieck space but 
there exists a sequence $(\mu_n)$ of Radon measures on $K$ such that
$\mu_n(A)\to0$ for every clopen set $A\subseteq K$ and $\int_Kfd\mu_n\not\to0$ for some $f\in C(K)$. Consequently, we get that Martin's axiom implies the existence of a Boolean algebra with the Grothendieck property but without the Nikodym property.
\end{abstract}


\maketitle

\section{Introduction}

A Banach space $X$ is \textit{Grothendieck} if every weak* convergent sequence in the dual space $X^*$ is also weakly convergent. The class of Grothendieck Banach spaces contains e.g. reflexive spaces, spaces $C(K)$ for basically disconnected compact Hausdorff spaces $K$ (Grothendieck \cite{Gro53}), the space $H^\infty$ of all bounded analytic functions on the unit disc (Bourgain \cite{Bou83}), von Neumann algebras (Pfitzner \cite{Pfi94}), etc. On the other hand, no non-reflexive separable Banach space is Grothendieck, thus no space $C(K)$ for $K$ metric compact is Grothendieck.

Transferring the definition of Grothendieck Banach spaces to the class of Boolean algebras, we say that a Boolean algebra $\aA$ has \textit{the Grothendieck property} if the Banach space $C(St(\aA))$ is a Grothendieck space, where $St(\aA)$ denotes the Stone space of $\aA$. Thus, by the aforementioned result of Grothendieck (\cite{Gro53}), every $\sigma$-complete Boolean algebra $\aA$ has the Grothendieck property as in this case $St(\aA)$ is basically disconnected (cf. \cite{vD80}). One can find further examples of Boolean algebras  with the Grothendieck property e.g. in \cite{Aiz88}, \cite{Das78}, \cite{Fre84vhs}, \cite{Hay81}, \cite{Hay01}, \cite{Mol81}, \cite{See68}, \cite{SZforext}. On the other hand, since for every countable Boolean algebra $\aA$ the Banach space $C(St(\aA))$ is separable, countable algebras do not have the property.

Recall that by the classical Riesz--Markov--Kakutani representation theorem, for every compact Hausdorff space $K$, the dual space $C(K)^*$ is isometrically isomorphic to the Banach space $M(K)$ of all signed Radon measures on $K$, endowed with the total variation norm. Furthermore, if $\aA$ is a Boolean algebra, $M(St(\aA))$ is isometrically isomorphic to the Banach space $\ba(\aA)$ of bounded signed finitely additive measures on $\aA$, also equipped with the total variation norm. Thus, $\aA$ has the Grothendieck property if and only if every  weak* convergent sequence in $\ba(\aA)$ is also weakly convergent.

A notion similar to the Grothendieck property is the Nikodym property. A Boolean algebra $\aA$ has \textit{the Nikodym property} if every pointwise bounded sequence $\seqn{\mu_n}$ in $\ba(\aA)$ is uniformly bounded, that is, $\sup_{n\io}|\mu_n(A)|<\infty$ for every $A\iA$ implies $\sup_{n\io}\|\mu_n\|<\infty$. Equivalently, $\aA$ has the Nikodym property if every sequence in $\ba(\aA)$ which is pointwise convergent on elements of $\aA$ is weak* convergent. Generalizing the result of Nikodym \cite{Nik33} for $\sigma$-additive measures on $\sigma$-fields, And\^{o} \cite{And61} proved that every $\sigma$-complete Boolean algebra has the Nikodym property. As in the case of the Grothendieck property, various weakenings of $\sigma$-completeness have also been established to imply the Nikodym property, see e.g. \cite{Das78}, \cite{Fre84vhs}, \cite{Hay81}, \cite{Mol81}, \cite{See68}. Similarly, it is easy to observe that countable Boolean algebras do not have the Nikodym property (see also \cite{MS24} and \cite{Zuc24} for more general results).

Schachermayer \cite{Sch82} found the first example of a Boolean algebra with the Nikodym property but without the Grothendieck property: the algebra $\mathcal{J}$ of Jordan measurable subsets of the unit interval $[0,1]$ (see also \cite{GW83} and \cite{Val13} for similar examples). Assuming Jensen's Diamond Principle $\Diamond$, the authors constructed in \cite{SZmingen} a minimally generated Boolean algebra with the Nikodym property---it is well-known known that minimally generated Boolean algebras cannot have the Grothendieck property (cf. \cite[Corollary 6.11]{KSZ}).

Finding an example of a Boolean algebra with the Grothendieck property but without the Nikodym property seems to be a much harder task as so far no such Boolean algebra has been obtained in ZFC, the standard system of axioms of set theory. 
Talagrand \cite{Tal84} constructed a Boolean algebra with the Grothendieck property but without the Nikodym property assuming the Continuum Hypothesis ($\mathsf{CH}$), that is, that the equality $\omega_1=\frakc$ holds, where $\omega_1$ denotes the first uncountable ordinal number and $\frakc$ is the cardinality of the real line $\mathbb{R}$. In this paper, following and exploiting the ideas of Talagrand \cite{Tal84}, we obtained the following strengthening of his result.

\begin{theorem}\label{theorem:main}
Assuming that $\frakp=\frakc$, there exists a Boolean algebra $\aA$ with the Grothendieck property but without the Nikodym property.
\end{theorem}

Recall that the \textit{pseudo-intersection number} $\frakp$ is defined as the minimal cardinality of a family $\fF$ of infinite subsets of the set $\N$ of natural numbers such that every non-empty finite subset of $\fF$ has infinite intersection and there is no infinite set $A\sub\N$ such that the difference $A\sm B$ is finite for every $B\in\fF$. Equivalently, $\frakp$ is the minimal character $\chi(X)$ of a countable topological space $X$ which is not Fr\'{e}chet--Urysohn. It holds $\omega_1\le\frakp\le\frakc$, yet none of the equalities $\omega_1=\frakp$ and $\frakp=\frakc$ can be proved in ZFC. We refer the reader to \cite{Bla10} and \cite{vD84} for detailed discussions on the number $\frakp$.

Since $\mathsf{CH}$ implies the equality $\frakp=\frakc$, Talagrand's result quickly follows from Theorem \ref{theorem:main}. The equality is also well-known to be a consequence of Martin's axiom $\mathsf{MA}$. In fact, by the theorem of Bell \cite{Bel81}, the equality $\frakp=\frakc$ is equivalent to a weak form of $\mathsf{MA}$, namely, to Martin's axiom for $\sigma$-centered posets, denoted $\mathsf{MA}(\sigma\text{-centered})$. Since $\mathsf{MA}$ implies $\mathsf{MA}(\sigma\text{-centered})$, we immediately get the next corollary. For information on Martin's axiom and its variants we refer the reader to \cite[Chapter 1]{BJ95} and \cite{Fre84}.

\begin{corollary}\label{cor:ma}
Assuming Martin's axiom, there exists a Boolean algebra $\aA$ with the Grothendieck property but without the Nikodym property.
\end{corollary}

No infinite Boolean algebra with the Grothendieck property can have cardinality strictly less than $\frakp$ (see \cite{SobKok}), therefore the assumption that $\frakp=\frakc$ implies that the Boolean algebra constructed in the proof of Theorem \ref{theorem:main} has cardinality $\frakc$. For a consistent example of a Boolean algebra of cardinality strictly less than $\frakc$ and with the Grothendieck property but without the Nikodym property, we refer the reader to the recent and independent work of G\l odkowski and Widz \cite{GW24}. Let us mention that their construction involves set-theoretic methods of forcing and was also  inspired by Talagrand's $\mathsf{CH}$ construction.

The pseudo-intersection number $\frakp$ has already appeared in the research connected with Grothendieck Banach spaces. Namely, Haydon, Levy, and Odell \cite{HLO87} proved, under the assumption that $\frakp=\frakc>\omega_1$, that every non-reflexive Grothendieck space has a quotient isomorphic to the Banach space $\elli$ of all bounded real-valued sequences. This contrasts another result of Talagrand \cite{Tal80} who, again assuming $\mathsf{CH}$ (so $\frakp=\frakc=\omega_1$), constructed a compact Hausdorff space $K$ such that the Banach space $C(K)$ is Grothendieck but it does not have any quotients isomorphic to $\elli$. Moreover, Haydon \cite{Hay81} constructed in ZFC a compact Hausdorff space $K$ such that $C(K)$ is a Grothendieck space but it does not contain any copy of $\elli$ (yet, it has a quotient isomorphic to $\elli$).

It is well-known that a Boolean algebra $\aA$ has the Nikodym property if and only if the linear dense subspace $S(St(\aA))$ of $C(St(\aA))$, spanned by the characteristic functions of clopen subsets of $St(\aA)$, is barrelled (cf. e.g. \cite[p. 10--11]{Sch82}). Thus, Theorem \ref{theorem:main} has the following geometric consequence.

\begin{corollary}\label{cor:ck_barrelled}
Assuming that $\frakp=\frakc$, there exists a totally disconnected compact Hausdorff space $K$ such that the Banach space $C(K)$ is a Grothendieck space but its subspace $S(K)$ is not barrelled. If moreover $\frakc>\omega_1$, then $C(K)$ has a quotient isomorphic to $\elli$.
\end{corollary}

Recall that by a result of R\"{a}biger \cite{Rab85}, a Banach space $C(K)$ is Grothendieck if and only if it does not have quotients isomorphic to the Banach space $c_0$ of all real-valued sequences converging to $0$. It follows that the space $C(K)$ from Corollary \ref{cor:ck_barrelled} does not have quotients isomorphic to $c_0$.


\section{Preliminaries, property (T), and auxiliary lemmas\label{sec:propt_lemmas}}

We first introduce some useful notations and prove auxiliary lemmas. We assume that $0\io$. As usual, $\Cantor$ denotes the Cantor space. For a compact Hausdorff space $K$, $Clopen(K)$ and $Bor(K)$ denote the Boolean algebras of all clopen and all Borel subsets of $K$, respectively, endowed with the standard set-theoretic operations. $C(K)$ denotes the Banach space of all continuous real-valued functions on $K$ endowed with the supremum norm and by $M(K)$ we denote the Banach space of all signed Radon measures on $K$ equipped with the total variation norm. Recall that by the classical Riesz--Markov--Kakutani representation theorem, $M(K)$ is isometrically isomorphic to the dual space $C(K)^*$. 

Let $\aA$ be a Boolean algebra. The zero element and the unit element of $\aA$ are denoted by $0_\aA$ and $1_\aA$, respectively. $St(\aA)$ denotes the Stone space of $\aA$, that is, the totally disconnected compact Hausdorff space consisting of all ultrafilters on $\aA$. By $\ba(\aA)$ we denote the Banach space of all bounded finitely additive measures on $\aA$ endowed with the total variation norm. When speaking about some measure $\mu$ on $\aA$, we tacitly assume that $\mu\in\ba(\aA)$. If $\mu\in\ba(\aA)$, then $|\mu|(\cdot)$ denotes the variation of $\mu$ and thus we have $\|\mu\|=|\mu|(1_\aA)$. Notice that $\ba(\aA)$ is isometrically isomorphic to $M(St(\aA))$.

For each $m\io$ and $i=0,1$ define the following clopen subset of $\Cantor$:
\[V_m^i=\{x\in\Cantor\colon\ x(m)=i\}.\]
For every $n\io$, let $\aA_n$ be the finite Boolean subalgebra of $Clopen(\Cantor)$ generated by the collection
\[\{V_m^i\colon\ 0\le m\le n, i=0,1\}.\]
Then, for every $n\io$, $\aA_n$ has $2^{n+1}$ atoms, and $Clopen(\Cantor)=\bigcup_{n\io}\aA_n$.

By $\lambda$ we mean the standard product measure on $\Cantor$, that is, $\lambda(V_m^i)=1/2$ for every $m\io$ and $i=0,1$. It follows that $\lambda(U)=2^{-(n+1)}$ for any atom $U$ of the algebra $\aA_n$, $n\io$. 

For each $n\io$ define the functions $\varphi_n\colon Bor(\Cantor)\to[0,\infty)$ and $\psi_n\colon Bor(\Cantor)\times Bor(\Cantor)\to[0,\infty)$ as follows:
\[\varphi_n(A)=|\lambda(A\cap V_n^0)-\lambda(A\cap V_n^1)|\]
and
\[\psi_n(A,B)=\big|(\lambda(A\cap V_n^0)-\lambda(A\cap V_n^1))-(\lambda(B\cap V_n^0)-\lambda(B\cap V_n^1))\big|,\]
for all $A,B\in Bor(\Cantor)$.

It is immediate that for every $A,B\in Bor(\Cantor)$ and $n\io$ we have
\[\psi_n(A,B)=\psi_n(B,A)\quad\text{and}\quad\psi_n(A,B)\le\varphi_n(A)+\varphi_n(B).\]
Also, for every $A\in Bor(\Cantor)$ and $n\io$, it holds $\psi_n(\emptyset,A)=\psi_n(A,\emptyset)=\varphi_n(A)$. The next three lemmas provide more inequalities for the functions $\psi_n$ and $\varphi_n$. Their proofs are left to the reader as they are routine.

\begin{lemma}\label{lemma:ineq}
For every $A,B,B',U\in Bor(\Cantor)$ and $m\io$ we have:
\begin{enumerate}[itemsep=1mm]
	\item $\lambda((A\cap B')\cap U)\le\lambda((A\cap B)\cap U)+\lambda(B'\triangle B)$,
	\item $\lambda(U\sm(A\cap B'))\le\lambda(U\sm(A\cap B))+\lambda(B'\triangle B)$,
	\item $\varphi_m((A\cap B')\cap U)\le\varphi_m((A\cap B)\cap U)+\psi_m\big((A\cap(B'\sm B))\cap U,\ (A\cap(B\sm B'))\cap U\big)$,
	\item $\psi_m\big((A\cap(B'\sm B))\cap U,\ (A\cap(B\sm B'))\cap U\big)\le\lambda(B'\triangle B)$,
	\item $\varphi_m((A\cap B')\cap U)\le\varphi_m((A\cap B)\cap U)+\lambda(B'\triangle B)$.
\end{enumerate}
\end{lemma}

\begin{lemma}\label{lemma:ineqpsi}
For every $W,W',Z,Z'\in Bor(\Cantor)$ such that $Z\cap Z'=\emptyset$ and $m\io$ we have:
\begin{enumerate}[itemsep=1mm]
	\item $\varphi_m(W\cap(Z\cup Z'))\le\varphi_m(W\cap Z)+\varphi_m(W\cap Z')$,
	\item $\psi_m\big(W\cap(Z\cup Z'),\ W'\cap(Z\cup Z')\big)\le\psi_m(W\cap Z,W'\cap Z)+\psi_m(W\cap Z',W'\cap Z')$.
\end{enumerate}
\end{lemma}

\begin{lemma}\label{lemma:ineqpsi2}
For every $W,W',Z,Z'\in Bor(\Cantor)$ and $m\io$ we have
\[\psi_m(W,Z)\le\psi_m(W',Z')+\lambda(W\triangle W')+\lambda(Z\triangle Z').\]
\end{lemma}

\subsection{Property (T)}

The property of Boolean subalgebras of $Clopen(\Cantor)$ introduced in the next definition is crucial for the proof of Theorem \ref{theorem:main} as it will imply among others that the algebra constructed therein will not have the Nikodym property (see Lemma \ref{lemma:no_nik}). Let us note that this property was first used by Talagrand \cite[page 166, point (1)]{Tal80}, however without giving it any name.

\begin{definition}\label{def:property_t} 
A subalgebra $\bB$ of $Bor(\Cantor)$ containing $Clopen(\Cantor)$ has \textit{property (T)} if for every finite sequence $A_1,\ldots,A_p\in\bB$ and $\eps>0$ there exists $n>0$ such that for every atom $U$ of $\aA_n$ and every $1\le i\le p$ the following two conditions hold:
\begin{enumerate}[(T.1)]
	\item either
	\[\lambda(A_i\cap U)\le\eps\lambda(U)/n,\]
	or
	\[\lambda(U\sm A_i)\le\eps\lambda(U)/n,\]
	\item for every $m>n$ we have
	\[\varphi_m(A_i\cap U)\le\eps\lambda(U)/m.\]
\end{enumerate}
\end{definition}

Intuitively speaking, elements of $\bB$ are ``well'' approximated by clopen subsets of $\Cantor$ with respect to the measure $\lambda$. Condition (T.1) might be understood as that the sets $A_i$ either almost ``omit'' the clopen $U$, or almost ``fill'' $U$, while condition (T.2) asserts that $A_i$ are ``symmetrical'' inside $U$ with respect to the measure $\lambda$ and the clopen subsets $V_m^0$ and $V_m^1$ for each $m>n$.


\begin{remark}\label{remark:def_prop_t}
Note that in the definition of property (T) we may actually require that $n$ is arbitrarily large, that is, if $\bB$ has property (T), then for every finite sequence $A_1,\ldots,A_p\in\bB$, $\eps>0$, and $N\io$ there exists $n\ge N$ such that for every atom $U$ of $\aA_n$ and every $1\le i\le p$ conditions (T.1) and (T.2) are satisfied. To see this, fix $A_1,\ldots,A_p\in\bB$, $\eps>0$, and $N\io$. Note that $\aA_N\sub Clopen(\Cantor)\sub\bB$ and use property (T) for the finite family $\fF=\{A_1,\ldots,A_p\}\cup\aA_N$ and $\eps'/2^{N+1}$, where $\eps'=\frac{1}{3}\min(1,\eps)$, to get $n>0$ as in Definition \ref{def:property_t}. If $n<N$, then for any atom $V$ of $\aA_N$ (so an element of $\fF$) there is an atom $U$ of $\aA_n$ such that $V\subsetneq U$. By condition (T.1), we then have either
\[2^{-(N+1)}=\lambda(V)=\lambda(V\cap U)\le(\eps'/2^{N+1})\cdot\lambda(U)/n<2^{-(N+1)}\cdot1/(2n)\le 2^{-(N+2)},\]
or 
\[2^{-(N+1)}\le\lambda(U\sm V)\le(\eps'/2^{N+1})\cdot\lambda(U)/n<2^{-(N+1)}\cdot1/(2n)\le 2^{-(N+2)},\]
so neither case is possible. It follows that $n\ge N$. Since $A_1,\ldots,A_p\in\fF$ and $\eps'/2^{N+1}<\eps$, conditions (T.1) and (T.2) are satisfied for $A_1,\ldots,A_p$ and $\eps$, too.

(In what follows, we will use the above remark frequently without mentioning it.)
\end{remark}

\begin{remark}\label{remark:def_prop_t_excl}
Note also that if for $\eps>0$ and $n>0$, $n\io$, we have $\eps/n<1/2$, then for any set $A\in Bor(\Cantor)$ and an atom $U$ of $\aA_n$ the alternative in condition (T.1) is exclusive, otherwise it would hold
\[\lambda(U)=\lambda(A\cap U)+\lambda(U\sm A)\le\eps\lambda(U)/n+\eps\lambda(U)/n<\lambda(U)/2+\lambda(U)/2=\lambda(U),\]
which is impossible.
\end{remark}

The following fact is obvious.

\begin{lemma}\label{lemma:union_prop_t}~
\begin{enumerate}
	\item The Boolean algebra $Clopen(\Cantor)$ has property (T).
	\item Let $\bB$ and $\bB'$ be subalgebras of $Bor(\Cantor)$ such that $Clopen(\Cantor)\sub\bB\sub\bB'$. If $\bB'$ has property (T), then $\bB$ has property (T).
	\item Let $\seq{\bB_\alpha}{\alpha<\beta}$ be an increasing chain of subalgebras of $Bor(\Cantor)$ such that $Clopen(\Cantor)\sub\bB_0$. If each $\bB_\alpha$ has property (T), then the union $\bigcup_{\alpha<\beta}\bB_\alpha$ has property (T).\noproof
\end{enumerate}
\end{lemma}

The next results imply that the whole $\sigma$-field $Bor(\Cantor)$ does not have property (T).

\begin{proposition}\label{prop:open_sets_prop_t}
Let $W$ be a non-empty clopen subset of $\Cantor$. Then, there are an open set $V_W\sub W$ and $N\io$ such that for every $n>N$, $n\io$, there is an atom $U$ of $\aA_n$ for which we have $\lambda(V_W\cap U)=\lambda(U)/2$ and $\lambda(U\sm V_W)=\lambda(U)/2$. Consequently, if a Boolean subalgebra $\bB$ of $Bor(\Cantor)$ contains $V_W$, then $\bB$ does not have property (T).
\end{proposition}
\begin{proof}
For any $k\io$ and $\sigma\in 2^{k+1}$ let $[\sigma]=\bigcap_{i=0}^{k}V_i^{\sigma(i)}$, so $[\sigma]$ is an atom of $\aA_k$.

Let $W$ be a non-empty clopen subset of $\Cantor$. There is $N\io$ and $\sigma\in 2^{N+1}$ such that $[\sigma]\sub W$. For every $k\io$ let $\tau_k\in 2^{(N+1)+k+2}$ be the sequence $\sigma1^k00$, that is, define
\[\tau_k(i)=\begin{cases}
\sigma(i),& 0\le i\le N,\\
1,& N< i\le N+k,\\
0,& i=N+k+1\ \text{ or }\ i=N+k+2,
\end{cases}\]
for every $0\le i\le N+k+2$. For each $k\io$, $\tau_k$ extends $\sigma$, so $[\tau_k]\sub[\sigma]\sub W$. $[\tau_k]$ is an atom of $\aA_{N+k+2}$. Set
\[V_W=\bigcup_{k\io}[\tau_k].\]
$V_W$ is an open subset of $\Cantor$ such that $V_W\sub[\sigma]\sub W$.

Fix $n>N$, $n\io$. Let $k=n-(N+1)$, so $(N+1)+k+2=n+2$, and hence $\tau_k\in 2^{n+2}$. Let $\rho\in2^{n+1}$ be such that $\rho(i)=\tau_k(i)$ for $0\le i\le n$ (so $\tau_k=\rho0$, hence $\tau_k$ extends $\rho$) and set $U=[\rho]$. Then, $U$ is an atom of $\aA_n$, $\lambda(U)=2^{-(n+1)}$, and $[\tau_k]\subsetneq U$. For every $m\neq k$ we also have $[\tau_m]\cap U=\emptyset$. Indeed, if $m>k$, then $\tau_m(N+k+1)=1$, whereas $\rho(N+k+1)=\tau_k(N+k+1)=0$, hence $[\tau_m]\cap U=[\tau_m]\cap[\rho]=\emptyset$. Similarly, if $m<k$, then $\rho(N+m+1)=\tau_k(N+m+1)=1$, whereas $\tau_m(N+m+1)=0$, hence again $[\tau_m]\cap U=[\tau_m]\cap[\rho]=\emptyset$. Consequently,
\[\lambda(V_W\cap U)=\lambda([\tau_k]\cap U)=\lambda([\tau_{k}])=2^{-(n+2)}=\lambda(U)/2\]
and so
\[\lambda(U\sm V_W)=\lambda(U)-\lambda(V_W\cap U)=\lambda(U)-\lambda(U)/2=\lambda(U)/2.\]

\medskip

The second part of the thesis follows by Remark \ref{remark:def_prop_t}.
\end{proof}

\begin{corollary}
No Boolean subalgebra of $Bor(\Cantor)$ with property (T) contains all open subsets of $\Cantor$.\noproof
\end{corollary}

The following lemma may seem a bit technical, but it will allow us to simplify the reasoning in a few places in the next section.

\begin{lemma}\label{lemma:complement}
Fix $n>0$, $n\io$, and $\eps>0$. Let $A,B,U\in Bor(\Cantor)$.
\begin{enumerate}[(a),itemsep=1mm]
	\item If
	\begin{itemize}
		\item either $\lambda(A\cap U)\le\eps\lambda(U)/2n$, or $\lambda(U\sm A)\le\eps\lambda(U)/2n$, and
		\item either $\lambda((A\cap B)\cap U)\le\eps\lambda(U)/2n$, or $\lambda(U\sm(A\cap B))\le\eps\lambda(U)/2n$,
	\end{itemize}
	then either $\lambda((A\sm B)\cap U)\le\eps\lambda(U)/n$, or $\lambda(U\sm(A\sm B))\le\eps\lambda(U)/n$.
	\item If
	$\varphi_m(A\cap U)\le\eps\lambda(U)/2m$ and $\varphi_m((A\cap B)\cap U)\le\eps\lambda(U)/2m$ for every $m>n$, 
	then $\varphi_m((A\sm B)\cap U)\le\eps\lambda(U)/m$ for every $m>n$.
\end{enumerate}
\end{lemma}
\begin{proof}
In order to prove (a) we need to consider several cases. If $\lambda(A\cap U)\le\eps\lambda(U)/2n$, then simply $\lambda((A\sm B)\cap U)\le\eps\lambda(U)/n$, so assume that $\lambda(U\sm A)\le\eps\lambda(U)/2n$. If $\lambda((A\cap B)\cap U)\le\eps\lambda(U)/2n$, then we have
\[\lambda(U\sm(A\sm B))=\lambda\big(U\sm(A\sm(A\cap B))\big)=\lambda\big(U\cap(A\cap(A\cap B)^c)^c\big)=\]
\[=\lambda\big(U\cap(A^c\cup(A\cap B))\big)=\lambda(U\sm A)+\lambda((A\cap B)\cap U)\le\]
\[\le\eps\lambda(U)/2n+\eps\lambda(U)/2n=\eps\lambda(U)/n.\]
Finally, assume that $\lambda(U\sm(A\cap B))\le\eps\lambda(U)/2n$. We get
\[\lambda(U)=\lambda(U\cap(A\cap B))+\lambda(U\sm(A\cap B))\le\lambda(U\cap(A\cap B))+\eps\lambda(U)/2n,\]
hence
\[\lambda(U\cap(A\cap B))\ge\lambda(U)-\eps\lambda(U)/2n.\]
It follows that
\[\lambda(U)\ge\lambda(U\cap A)=\lambda(U\cap(A\cap B))+\lambda(U\cap(A\sm B))\ge\]
\[\ge\lambda(U)-\eps\lambda(U)/2n+\lambda(U\cap(A\sm B)),\]
whence, after rearranging, we get
\[\lambda((A\sm B)\cap U)\le\eps\lambda(U)/2n<\eps\lambda(U)/n,\]
which finishes the proof of (a).

We now prove that (b) holds as well. Let $m>n$. We have
\[\varphi_m((A\sm B)\cap U)=\big|\lambda(A\cap B^c\cap U\cap V_m^0)-\lambda(A\cap B^c\cap U\cap V_m^1)\big|=\]
\[=\big|\lambda(A\cap U\cap V_m^0)-\lambda(A\cap B\cap U\cap V_m^0)-\lambda(A\cap U\cap V_m^1)+\lambda(A\cap B\cap U\cap V_m^1)\big|\le\]
\[\le\big|\lambda(A\cap U\cap V_m^0)-\lambda(A\cap U\cap V_m^1)\big|+\big|\lambda(A\cap B\cap U\cap V_m^0)-\lambda(A\cap B\cap U\cap V_m^1)\big|=\]
\[=\varphi_m(A\cap U)+\varphi_m((A\cap B)\cap U)\le\eps\lambda(U)/2m+\eps\lambda(U)/2m=\eps\lambda(U)/m,\]
so (b) holds, too.
\end{proof}

The next result expresses an important feature of property (T), namely, that it induces the existence of nice homomorphisms from finite subalgebras of Boolean algebras with property (T) into algebras of the form $\aA_n$.

\begin{lemma}\label{lemma:homomorphism}
Let $n>0$, $n\io$, and $0<\eps<1/4$. Let $\eE$ be a finite subalgebra of $Bor(\Cantor)$ such that for every $A\in\eE$ and every atom $V$ of $\aA_n$ either $\lambda(A\cap V)\le\eps\lambda(V)/n$, or $\lambda(V\sm A)\le\eps\lambda(V)/n$. Let $\tau\colon\eE\to\aA_n$ be a function defined for every $A\in\eE$ as follows:
\[\tau(A)=\bigcup\big\{V\colon\ V\text{ is an atom of }\aA_n\text{ such that }\lambda(V\sm A)\le\eps\lambda(V)/n\big\}.\]
Then, $\tau$ is a Boolean homomorphism such that for every $A\in\eE$ we have
\[\lambda(A\triangle\tau(A))\le\eps/n.\]
\end{lemma}
\begin{proof}
Notice first that by Remark \ref{remark:def_prop_t_excl} for any $A\in\eE$ and any atom $V$ of $\aA_n$ it is impossible that both inequalities $\lambda(A\cap V)\le\eps\lambda(V)/n$ and $\lambda(V\sm A)\le\eps\lambda(V)/n$ hold simultaneously.

(1) Obviously, for the zero and unit elements we have $\tau(0_\eE)=\tau(\emptyset)=\emptyset=0_{\aA_n}$ and $\tau(1_\eE)=\tau(\Cantor)=\Cantor=1_{\aA_n}$. 

\medskip

(2) Let $A\in\eE$. By the definition of $\tau$, we have
\[\tau(A^c)=\bigcup\big\{V\colon\ V\text{ is an atom of }\aA_n\text{ such that }\lambda(V\sm A^c)\le\eps\lambda(V)/n\big\}.\]
For an atom $V$ of $\aA_n$ note that $V\sm A^c=A\cap V$, so if $V\sub\tau(A^c)$, then, by the remark at the beginning of the proof, we have $V\sub\tau(A)^c$. It follows that $\tau(A^c)\sub\tau(A)^c$.
We also have
\[\tau(A)^c=\bigcup\big\{V\colon\ V\text{ is an atom of }\aA_n\text{ such that }\lambda(A\cap V)\le\eps\lambda(V)/n\big\}.\]
Simiarly as above, for an atom $V$ of $\aA_n$ it holds $A\cap V=V\sm A^c$, so if $V\in\tau(A)^c$, then $V\in\tau(A^c)$. It follows that $\tau(A)^c\sub\tau(A^c)$. Hence, $\tau(A^c)=\tau(A)^c$.

\medskip

(3) Let $A,B\in\eE$. If for an atom $V$ of $\aA_n$ we have $\lambda(V\sm A)\le\eps\lambda(V)/n$ or $\lambda(V\sm B)\le\eps\lambda(V)/n$, then trivially $\lambda(V\sm (A\cup B))\le\eps\lambda(V)/n$. It follows that $\tau(A)\cup\tau(B)\sub\tau(A\cup B)$. We thus need only to show that $\tau(A\cup B)\sub\tau(A)\cup\tau(B)$.

Let $V$ be an atom of $\aA_n$ such that $\lambda(V\sm(A\cup B))\le\eps\lambda(V)/n$, that is, $V\sub\tau(A\cup B)$. Assume that $V\not\subseteq\tau(A)\cup\tau(B)$. It follows that $V\cap(\tau(A)\cup\tau(B))=\emptyset$, hence $\lambda(V\sm A)>\eps\lambda(V)/n$ and $\lambda(V\sm B)>\eps\lambda(V)/n$. Consequently, $\lambda(A\cap V)\le\eps\lambda(V)/n$ and $\lambda(B\cap V)\le\eps\lambda(V)/n$. We have
\[\lambda((A\cup B)\cap V)=\lambda((A\cap V)\cup(B\cap V))\le\lambda(A\cap V)+\lambda(B\cap V)\le2\eps\lambda(V)/n,\]
so, similarly as at the beginning of the proof, we get
\[\lambda(V)=\lambda((A\cup B)\cap V)+\lambda(V\sm(A\cup B))\le2\eps\lambda(V)/n+\eps\lambda(V)/n=\]
\[=3\eps\lambda(V)/n<(3/4)\lambda(V)<\lambda(V),\]
which is a contradiction. It thus follows that $V\sub\tau(A)\cup\tau(B)$.

Having proved that $\tau(A\cup B)\sub\tau(A)\cup\tau(B)$, we get that $\tau(A\cup B)=\tau(A)\cup\tau(B)$.

\medskip

Paragraphs (1)--(3) prove that $\tau$ is a homomorphism of the Boolean algebras. To prove the last statement, fix $A\in\eE$. We have
\[\lambda(A\triangle\tau(A))=\sum_{V\text{ -- atom of }\aA_n}\lambda((A\triangle\tau(A))\cap V)=\]
\[=\sum_{V\text{ -- atom of }\aA_n}\big(\lambda((A\sm\tau(A))\cap V)+\lambda((\tau(A)\sm A)\cap V)\big)=\]
\[=\sum_{V\text{ -- atom of }\aA_n}\big(\lambda(A\cap\tau(A)^c\cap V)+\lambda(\tau(A)\cap A^c\cap V)\big)=\]
\[=\sum_{\substack{V\text{ -- atom of }\aA_n\\\lambda(V\sm A)\le\eps\lambda(V)/n}}\big(\lambda(A\cap\tau(A)^c\cap V)+\lambda(\tau(A)\cap A^c\cap V)\big)+\]
\[+\sum_{\substack{V\text{ -- atom of }\aA_n\\\lambda(V\cap A)\le\eps\lambda(V)/n}}\big(\lambda(A\cap\tau(A)^c\cap V)+\lambda(\tau(A)\cap A^c\cap V)\big).\]
If for an atom $V$ of $\aA_n$ we have $\lambda(V\sm A)\le\eps\lambda(V)/n$, then $V\sub\tau(A)$ and hence $\lambda(A\cap\tau(A)^c\cap V)=\lambda(\emptyset)=0$. Similarly, if $\lambda(V\sm A^c)=\lambda(V\cap A)\le\eps\lambda(V)/n$, then $V\sub\tau(A^c)=\tau(A)^c$, hence $\lambda(\tau(A)\cap A^c\cap V)=\lambda(\emptyset)=0$. Thus, we get
\[\lambda(A\triangle\tau(A))=\sum_{\substack{V\text{ -- atom of }\aA_n\\\lambda(V\sm A)\le\eps\lambda(V)/n}}\lambda(\tau(A)\cap A^c\cap V)+\sum_{\substack{V\text{ -- atom of }\aA_n\\\lambda(V\cap A)\le\eps\lambda(V)/n}}\lambda(A\cap\tau(A)^c\cap V)\le\]
\[\le\sum_{\substack{V\text{ -- atom of }\aA_n\\\lambda(V\sm A)\le\eps\lambda(V)/n}}\lambda(V\sm A)+\sum_{\substack{V\text{ -- atom of }\aA_n\\\lambda(V\cap A)\le\eps\lambda(V)/n}}\lambda(V\cap A)\]
\[\le\big(\lambda(\tau(A))/2^{-(n+1)}\big)\cdot\eps2^{-(n+1)}/n+\big(2^{n+1}-\lambda(\tau(A))/2^{-(n+1)}\big)\cdot\eps2^{-(n+1)}/n=\eps/n,\]
where the last inequality follows from the fact that the number of atoms of $\aA_n$ contained in $\tau(A)$ is equal to $\lambda(\tau(A))/2^{-(n+1)}$.
\end{proof}

The following lemma also presents an important attribute of Boolean algebras having property (T): they do not have the Nikodym property.

\begin{lemma}\label{lemma:no_nik}
Let $\bB$ be a Boolean subalgebra of $Bor(\Cantor)$ containing $Clopen(\Cantor)$. If $\bB$ satisfies property (T), then $\bB$ does not have the Nikodym property.
\end{lemma}
\begin{proof}
For each $n\io$ and $A\in\bB$ set
\[\mu_n(A)=n(\lambda(A\cap V_n^0)-\lambda(A\cap V_n^1));\]
then, $\mu_n\in\ba(\bB)$ and $\|\mu_n\|=n$ (recall that $V_n^0$ and $V_n^1$ are disjoint and belong to $Clopen(\Cantor)$ and hence to $\bB$).

We show that $\seqn{\mu_n}$ is pointwise convergent to $0$ (and so that it is pointwise bounded). Fix $A\in\bB$ and $\eps>0$. By property (T) of $\bB$, there is $n\io$ such that for all atoms $U\in\aA_n$ and $m>n$ we have
\[\varphi_m(A\cap U)\le\eps\lambda(U)/m.\]
Fix $m>n$. Since $A=\bigcup\{A\cap U\colon\ U\text{ -- atom of }\aA_n\}$, we get that
\[\varphi_m(A)=|\lambda(A\cap V_m^0)-\lambda(A\cap V_m^1)|=\big|\sum_{U\text{ -- atom of }\aA_n}(\lambda(A\cap V_m^0\cap U)-\lambda(A\cap V_m^1\cap U))\big|\le\]
\[\le\sum_{U\text{ -- atom of }\aA_n}\big|\lambda(A\cap V_m^0\cap U)-\lambda(A\cap V_m^1\cap U)\big|=\sum_{U\text{ -- atom of }\aA_n}\varphi_m(A\cap U)\le\]
\[\le\eps/m\cdot\sum_{U\text{ -- atom of }\aA_n}\lambda(U)=\eps/m.\]
Since $|\mu_m(A)|=m\cdot\varphi_m(A)$, we have
\[|\mu_m(A)|=m\cdot\varphi_m(A)\le m\cdot\eps/m=\eps,\]
which implies that $\mu_n(A)\to0$ as $n\to\infty$. 

It follows that $\seqn{\mu_n}$ witnesses that $\bB$ does not have the Nikodym property.
\end{proof}

\noindent (Note that in the above proof we only use condition (T.2) of property (T).)

\subsection{Characterization of the Grothendieck property}

We finish the section with the following characterization of the Grothendieck property of Boolean algebras. It is already part of the folklore (cf. e.g. \cite{Tal84}, \cite[Lemma 2.2]{KS13}), but we provide its proof for the sake of completeness and because of the apparent lack of a direct proof in the literature. Recall that a collection $\seq{A_i}{i\in I}$ of elements of a Boolean algebra $\aA$ is \textit{an antichain} if $A_i\wedge A_j=0_{\aA}$ for every $i\neq j\in I$. 

\begin{lemma}\label{lemma:kpr}
A Boolean algebra $\aA$ has the Grothendieck property if and only if $\aA$ satisfies the following condition:
\begin{itemize}[\textup{($\dagger$)}]
\item for every $\eps\in(0,1)$ and every sequence $\seqn{\nu_n}$ in $\ba(\aA)$ with $\|\nu_n\|=1$ for each $n\io$, if there exists an antichain $\seqn{H_n}$ in $\aA$ such that $|\nu_n|(H_n)\ge1-\eps$ for each $n\io$, then there also exists an element $B\in\aA$ such that the limit $\lim_{n\to\infty}\nu_n(B)$ does not exist.
\end{itemize}%
%
\end{lemma}
\begin{proof}
Assume first that $\aA$ has property ($\dagger$) but it does not have the Grothendieck property, that is, that $C(St(\aA))$ is not a Grothendieck space. There is a weak* null sequence $\seqn{\mu_n}$ in $M(St(\aA))$ such that none of its subsequences is weakly null. Without loss of generality we may assume that $\|\mu_n\|=1$ for every $n\io$. By Kadec--Pe\l czy\'{n}ski--Rosenthal's Subsequential Splitting Lemma (cf. \cite[Lemma 5.2.8]{AK06}), there are $r\ge0$, a subsequence $\seqk{\mu_{n_k}}$, and sequences $\seqk{\lambda_k}$ and $\seqk{\theta_k}$ in $M(St(\aA))$, such that:
\begin{itemize}
	\item $\mu_{n_k}=\lambda_k+\theta_k$ for every $k\io$,
	\item $\seqk{\lambda_k}$ is supported by disjoint Borel subsets of $St(\aA)$,
	\item $\seqk{\theta_k}$ is weakly convergent,
	\item $\|\lambda_k\|=r$ for every $k\io$.
\end{itemize}
If $r=0$, then $\mu_{n_k}=\theta_k$ for each $k\io$, and so $\seqk{\mu_{n_k}}$ is weakly null, a contradiction. Hence, $r>0$. For each $k\io$ set $\rho_k=\lambda_k/r$. It follows that $\seqk{\rho_k}$ is supported by disjoint Borel subsets of $St(\aA)$ and $\|\rho_k\|=1$ for every $k\io$. By \cite[Lemma 1]{Tal80} and the regularity of $\rho_k$'s, there are a subsequence $\seql{\rho_{k_l}}$ and an antichain $\seql{H_l}$ in $\aA$ such that $|\rho_{k_l}|([H_l])\ge1-\eps$, where, for $A\iA$, $[A]$ denotes the clopen subset of $St(\aA)$ corresponding via the Stone duality to $A$. For each $l\io$ let $\nu_l\in\ba(\aA)$ be such that $\nu_l(A)=\rho_{k_l}([A])$ for every $A\iA$. Then, $\|\nu_l\|=\|\rho_{k_l}\|$ and $|\nu_l|(H_l)\ge1-\eps$ for every $l\io$. By ($\dagger$), there is $B\in\aA$ such that $\lim_{l\to\infty}\nu_l(B)$ does not exist. Consequently, $\lim_{l\to\infty}\rho_{k_l}([B])$ does not exist. But this is a contradiction, since $\seqk{\lambda_k}$ is weak* convergent (being a difference of a weak* null sequence and a weakly convergent sequence) and hence so is also $\seql{\rho_{k_l}}$, which in particular means that $\seql{\rho_{k_l}}$ is convergent on clopen subsets of $St(\aA)$. It follows that $\aA$ has the Grothendieck property.

Assume now that $\aA$ has the Grothendieck property, we will show that it has property ($\dagger$) as well. Let $\eps\in(0,1)$ and $\seqn{\nu_n}$ be a sequence in $\ba(\aA)$ such that $\|\nu_n\|=1$ for every $n\io$ and there exists an antichain $\seqn{H_n}$ in $\aA$ for which we have $|\nu_n|(H_n)\ge1-\eps$ for every $n\io$. We will work in $M(St(\aA))$, that is, for each $n\io$ let $\rho_n\in M(St(\aA))$ be a unique Radon measure such that $\rho_n([A])=\nu_n(A)$ for every $A\in\aA$ (where $[\cdot]$ has the same meaning as above). Let $\delta=1-\eps$. By the regularity of each $\rho_n$, there exists an antichain $\seqn{G_n}$ in $\aA$ such that $G_n\le H_n$ and $|\rho_n([G_n])|>\delta/3$ for every $n\io$. By Rosenthal's lemma (\cite[page 82]{Die84}, cf. also \cite{SobRos}, \cite{KMC21}), there is a sequence $\seqk{n_k\io}$ such that for every $k\io$ we have
\[\tag{$*$}\sum_{\substack{l\io\\l\neq k}}|\rho_{n_k}|([G_{n_l}])<\delta/9.\]
Set $B=\bigcup_{l\io}[G_{n_{2l}}]$. It follows that $B$ is a Borel (open) subset of $St(\aA)$, so the functional $\varphi$ on $M(St(\aA))$ defined as follows:
\[M(St(\aA))\ni\rho\xmapsto{\,\ \varphi\ \,}\int_{St(\aA)}\chi_Bd\rho,\]
where $\chi_B$ denotes the characteristic function of $B$ in $St(\aA)$, is in $M(St(\aA))^*$. By ($*$), for each $k\io$ we have
\[|\varphi(\rho_{n_{2k}})|=\big|\int_{St(\aA)}\chi_Bd\rho_{n_{2k}}\big|\ge\big|\int_{St(\aA)}\chi_{[G_{n_{2k}}]}d\rho_{n_{2k}}\big|-\big|\int_{St(\aA)}\chi_{B\sm[G_{n_{2k}}]}d\rho_{n_{2k}}\big|\ge\]
\[\ge\big|\int_{St(\aA)}\chi_{[G_{n_{2k}}]}d\rho_{n_{2k}}\big|-\int_{St(\aA)}\chi_{B\sm[G_{n_{2k}}]}d|\rho_{n_{2k}}|=|\rho_{n_{2k}}([G_{n_{2k}}])|-\sum_{\substack{l\io\\l\neq k}}|\rho_{n_{2k}}|([G_{n_{2l}}])>\]
\[>\delta/3-\delta/9=2\delta/9,\]
as well as
\[|\varphi(\rho_{n_{2k+1}})|=\big|\int_{St(\aA)}\chi_Bd\rho_{n_{2k+1}}\big|\le\int_{St(\aA)}\chi_Bd|\rho_{n_{2k+1}}|=\sum_{l\io}|\rho_{n_{2k+1}}|(G_{n_{2l}})<\delta/9.\]
Consequently,
\[\limsup_{n\to\infty}|\varphi(\rho_n)|\ge2\delta/9>\delta/9\ge\liminf_{n\to\infty}|\varphi(\rho_n)|,\]
so $\seqn{\rho_n}$ is not weakly convergent. Since $C(St(\aA))$ is a Grothendieck space and the spaces $C(St(\aA))^*$ and $M(St(\aA))$ are isomorphic, $\seqn{\rho_n}$ is not weak* convergent either. As $\sup_{n\io}\|\rho_n\|=1<\infty$, by the Stone--Weierstrass theorem, there is an element $B\in\aA$ such that the limit $\lim_{n\to\infty}\rho_n([B])$ does not exist. Consequently, the limit $\lim_{n\to\infty}\nu_n(B)$ does not exist. It follows that $\aA$ has property ($\dagger$).
\end{proof}

We will use the above lemma in the proof of Theorem \ref{theorem:main} in Section \ref{sec:proof} with $\eps=0.05$ (so $1-\eps=0.95$).

\section{Extending algebras with property (T)\label{sec:ext_prop_t}}

%
%
We will work with the following setting.  Let $\cC\sub\bB$ be Boolean subalgebras of $Bor(\Cantor)$ such that the following conditions are satisfied:
\begin{itemize}
	\item $\cC$ is countable,
	\item $\cC$ contains $Clopen(\Cantor)$,
	\item and $\bB$ satisfies property (T). 
\end{itemize}
Let $\seqn{\nu_n}$ be a sequence of measures on $\cC$ and $\seqn{H_n}$ an antichain in $\cC$ such that $\|\nu_n\|=1$ and $|\nu_n|(H_n)\ge0.95$ for every $n\io$.

Assume also that there is a non-negative measure $\theta$ on $\cC$ such that $\lim_{n\to\infty}|\nu_n|(A)=\theta(A)$ for every $A\in\cC$. By the Lebesgue Decomposition Theorem for finitely additive measures on fields (see \cite[Theorem 6.2.4]{RR83}), there exist two finitely additive non-negative measures $\theta_1$ and $\theta_2$ on $\cC$ (i.e. $\theta_1,\theta_2\in\ba(\cC)$ and $\theta_1,\theta_2\ge0$) such that the following hold:
\begin{itemize}
	\item $\theta=\theta_1+\theta_2$,
	\item $\theta_1$ is \textit{absolutely continuous} with respect to $\lambda$ (i.e. for every $\eps>0$ there is $\delta>0$ such that for every $A\in\cC$ if $\lambda(A)<\delta$, then $\theta_1(A)<\eps$), and 
	\item $\theta_2$ is \textit{orthogonal} to $\lambda$ (i.e., for every $\eps>0$ there is $X\in\cC$ such that $\lambda(X)<\eps$ and $\theta_2(X^c)<\eps$).
\end{itemize}


\medskip

%
%

We start by recalling the following result from Talagrand's article \cite{Tal84}. Since this lemma is crucial for our construction (as it was in the case of \cite{Tal84}), we include its technical proof in the appendix at the end of the paper.

\begin{lemma}[{\cite[{\textit{Lemme}}]{Tal84}}]\hspace{-2mm}
\label{lemma:combinatorial}
Let $t\io$, $t>0$, and $0<\eta<2^{-t-11}$. There exists $n_0\ge t$ such that for every
\begin{itemize}
	\item $n\ge n_0$, 
	\item sets $P,R,Z\in\aA_n$ with $\max(\lambda(P),\lambda(R),\lambda(Z))\le2^{-8}\cdot\eta^2$, and
	\item set $Q\in\aA_n$ for which there exists an atom $G$ of $\aA_t$ with $\lambda(Q\cap G)\ge0.95\lambda(G)$,
\end{itemize}
there exists a subset $M\in\aA_n$ of $Q$ such that $\lambda(M)\le\eta$, $M\cap(P\cup R\cup Z)=\emptyset$, and
\[\psi_m(M\cup(Z\cap Q),\ R\cap Q)\le\eta/m\]
for every $t<m\le n$.
\end{lemma}
\begin{proof}
See Appendix.
\end{proof}

\begin{lemma}
\label{lemma:troisieme}
Let $\dD_0$ be a finite collection of pairwise disjoint elements of $\bB$ such that $|\dD_0|\ge2$. Let $t\io$, $t>0$, be such that for every element $C\in\dD_0$ there is an atom $U_C$ of $\aA_t$ for which we have $\lambda(U_C\cap C)>0.99\lambda(U_C)$. Let $0<\eta<2^{-t-11}$ and $\zeta=2^{-10}\eta^2/|\dD_0|^2$. Let $X\in\cC$ be such that $\lambda(X)<\zeta$. 

Then, for every $W,Y,H\in\cC$ with $W,Y\sub H$ and  $\lambda(H)<\zeta$, there is $A\in\cC$ such that the following conditions hold:
\begin{enumerate}[({A}.1)]
	\item $A\cap (X\cup H)=\emptyset$,
	\item $\lambda(A)\le\eta$,
	\item for every $C\in\dD_0$ and $m>t$ we have
	\[\psi_m\big((A\cup Y)\cap C,\ W\cap C\big)\le\eta/m.\]
\end{enumerate}
\end{lemma}
\begin{proof}
Let $n_0\ge t$ be as in Lemma \ref{lemma:combinatorial} (for $\eta/|\dD_0|$ in place of $\eta$). Fix $W,Y,H\in\cC$ with $W,Y\sub H$ and $\lambda(H)<\zeta$. 

Let $\eE$ denote the finite Boolean subalgebra of $\bB$ generated by the family \[\dD_0\cup\aA_t\cup\{X,W,Y,H\}.\]
Let
\[\delta=\eta^2/(30|\dD_0|^2).\]
Note that $\delta\le\eta^2/30\le\eta/30$. By property (T) of $\bB$, there is $n\ge\max(n_0,30)$ such that for every $E\in\eE$ and every atom $V$ of $\aA_n$ the following hold:
\begin{enumerate}[{(C.1)}]
	\item either
	\[\lambda(E\cap V)\le\delta\lambda(V)/n,\]
	or
	\[\lambda(V\sm E)\le\delta\lambda(V)/n,\]
	\item for every $m>n$ we have
	\[\varphi_m(E\cap V)\le\delta\lambda(V)/m.\]
\end{enumerate}
Let $\tau\colon\eE\to\aA_n$ be a function defined for every $E\in\eE$ as follows:
\[\tau(E)=\bigcup\big\{V\colon\ V\text{ is an atom of }\aA_n\text{ such that }\lambda(V\sm E)\le\delta\lambda(V)/n\big\}.\]
Lemma \ref{lemma:homomorphism} (with $\eps=\delta$) implies that $\tau$ is a homomorphism of Boolean algebras such that
\[\tag{C.3}\lambda(E\triangle\tau(E))\le\delta/n\le\eta/30n\]
for every $E\in\eE$.

\medskip

For every $C\in\dD_0$ set also
\[Q_C=\tau(C).\]
We claim that $\lambda(Q_C\cap U_C)\ge0.95\lambda(U_C)$. To see this, note that 
\[\delta/n\le\eta/30\le2^{-t-11}/30<0.001\lambda(U_C),\]
so, by (C.3) and the fact that $\tau(U_C)=U_C$ (as $U_C$ is an element of $\aA_t\sub\aA_n$), we have
\[\lambda(Q_C\cap U_C)=\lambda(\tau(C)\cap U_C)=\lambda(\tau(C\cap U_C))\ge\lambda(C\cap U_C)-\lambda((C\cap U_C)\sm\tau(C\cap U_C))\ge\]
\[\ge0.99\lambda(U_C)-\delta/n>0.99\lambda(U_C)-0.001\lambda(U_C)>0.95\lambda(U_C).\]

Set $Z=\tau(Y)$, $R=\tau(W)$, and $P=\tau(X\cup H)$. Since $n\ge30$, we have
\[\lambda(Z)=\lambda(\tau(Y))=\lambda(\tau(Y)\cap Y)+\lambda(\tau(Y)\sm Y)\le\lambda(Y)+\lambda(Y\triangle\tau(Y))<\]
\[<\zeta+\delta/n=2^{-10}\eta^2/|\dD_0|^2+\eta^2/(30|\dD_0|^2\cdot n)\le2^{-8}(\eta/|\dD_0|)^2,\]
so $\lambda(Z)\le2^{-8}(\eta/|\dD_0|)^2$. Similarly, $\max(\lambda(W),\lambda(P))\le2^{-8}(\eta/|\dD_0|)^2$.

\medskip

By the property of $n_0$ given by Lemma \ref{lemma:combinatorial}, for each $C\in\dD_0$ there is a subset $M_C\sub Q_C$ such that $M_C\in\aA_n$, $\lambda(M_C)\le\eta/|\dD_0|$, $M_C\cap(P\cup R\cup Z)=\emptyset$, and
\[\tag{C.4}\psi_m\big(M_C\cup(Z\cap Q_C),\ R\cap Q_C\big)\le\eta/(|\dD_0|\cdot m)\]
for every $t<m\le n$. Since $\tau$ is a homomorphism, $M_C\cap M_{C'}=\emptyset$ for every $C\neq C'\in\dD_0$. Put
\[A=\bigcup_{C\in\dD_0}M_C\sm(X\cup H).\]
Obviously, $A\cap(X\cup H)=\emptyset$ and $\lambda(A)\le\eta$, so conditions (A.1) and (A.2) hold. We only need to show that condition (A.3) holds, too.

\medskip

Let $m>t$. Assume first that $m>n$. Note that the set $A$ is of the form $(X\cup H)^c\cap M$, where $(X\cup H)^c\in\eE$ and $M=\bigcup_{C\in\dD_0}M_C\in\aA_n$ (so $M$ is a union of atoms of $\aA_n$). Also, $(X\cup H)^c\cap C\in\eE$ and $Y\cap C\in\eE$ for every $C\in\dD_0$. Consequently, by Lemma \ref{lemma:ineqpsi}.(1) (with $W=(X\cup H)^c\cap C$, and $Z$ and $Z'$ ranging over all atoms of $\aA_n$ contained in $M$) and condition (C.2), for every $C\in\dD_0$ we have
\[\varphi_m(M\cap(X\cup H)^c\cap C)\le\sum_{V\text{ -- atom of }\aA_n}\delta\lambda(V)/m=\delta/m,\]
and, by a similar argument, also
\[\varphi_m(Y\cap C)\le\delta/m.\]
Thus, since $A\cap Y=\emptyset$ (as $Y\sub H$), again by Lemma \ref{lemma:ineqpsi}.(1) (with $W=C$, $Z=A$, and $Z'=Y$), we have
\[\varphi_m((A\cup Y)\cap C)\le\varphi_m(A\cap C)+\varphi_m(Y\cap C)=\]
\[=\varphi_m(M\cap(X\cup H)^c\cap C)+\varphi_m(Y\cap C)\le\]
\[\le\delta/m+\delta/m=2\delta/m.\]
Similarly, since $W\cap C\in\eE$ for every $C\in\dD_0$, by Lemma \ref{lemma:ineqpsi}.(1) (with $W\cap C$ in place of $W$, and $Z$ and $Z'$ ranging over all atoms of $\aA_n$) and condition (C.2), for every $C\in\dD_0$ we have
\[\varphi_m(W\cap C)\le\delta/m.\] 
The last two inequalities together yield
\[\psi_m((A\cup Y)\cap C,\ W\cap C)\le\varphi_m((A\cup Y)\cap C)+\varphi_m(W\cap C)\le\]
\[\le2\delta/m+\delta/m=3\delta/m=\eta^2/(10|\dD_0|^2\cdot m)\le\eta/m,\]
so condition (A.3) is satisfied for $m>n$.

\medskip

Assume now that $t<m\le n$. Fix $C\in\dD_0$. We need to do a series of estimations. Note first that
\[\tag{C.5}M_C\cup(Z\cap Q_C)=(M_C\cup Z)\cap Q_C.\]
Hence, it holds
\[((A\cup Y)\cap C)\sm(M_C\cup(Z\cap Q_C))=((A\cup Y)\cap C)\sm((M_C\cup Z)\cap Q_C).\]
We will estimate the value of $\lambda\big(((A\cup Y)\cap C)\sm((M_C\cup Z)\cap Q_C)\big)$. We have
\[((A\cup Y)\cap C)\sm((M_C\cup Z)\cap Q_C)=(A\cup Y)\cap C\cap((M_C\cup Z)^c\cup Q_C^c)=\]
\[=\big((A\cup Y)\cap C\cap(M_C\cup Z)^c\big)\cup\big((A\cup Y)\cap C\cap Q_C^c\big)=\]
\[=\big((A\cup Y)\cap C\cap M_C^c\cap Z^c\big)\cup\big((A\cup Y)\cap C\cap Q_C^c\big)=\]
\[=\big(A\cap C\cap M_C^c\cap Z^c\big)\cup\big(Y\cap C\cap M_C^c\cap Z^c\big)\cup\big((A\cup Y)\cap C\cap Q_C^c\big)\sub\]
\[\sub(A\cap C\cap M_C^c)\cup(Y\cap Z^c)\cup(C\cap Q_C^c).\]
Since $A\cap M_C^c\sub Q_C^c$, for $A\cap C\cap M_C^c$ we have
\[A\cap C\cap M_C^c\sub C\cap Q_C^c=C\sm\tau(C),\]
so, by (C.3), we have
\[\tag{D.1}\lambda\big((A\cap C\cap M_C^c)\cup(C\cap Q_C^c)\big)=\lambda(C\cap Q_C^c)\le\eta/30n.\]
Since $Z=\tau(Y)$, condition (C.3) also implies that
\[\tag{D.2}\lambda(Y\cap Z^c)\le\eta/30n.\]
Inequalities (D.1) and (D.2) yield
\[\tag{E.1}\lambda\big(((A\cup Y)\cap C)\sm((M_C\cup Z)\cap Q_C)\big)\le2\eta/30n.\]

Second, we will estimate $\lambda\big(((M_C\cup Z)\cap Q_C)\sm((A\cup Y)\cap C)\big)$. We have
\[((M_C\cup Z)\cap Q_C)\sm((A\cup Y)\cap C)=(M_C\cup Z)\cap Q_C\cap((A\cup Y)^c\cup C^c)=\]
\[=\big((M_C\cup Z)\cap Q_C\cap(A\cup Y)^c\big)\cup\big((M_C\cup Z)\cap Q_C\cap C^c\big)\sub\]
\[\sub\big((M_C\cup Z)\cap Q_C\cap(A\cup Y)^c\big)\cup(Q_C\cap C^c).\]
Since $Q_C\cap C^c=\tau(C)\sm C$, by (C.3), we have
\[\tag{D.3}\lambda(Q_C\cap C^c)\le\eta/30n.\]
It holds
\[(M_C\cup Z)\cap Q_C\cap(A\cup Y)^c=(M_C\cup Z)\cap Q_C\cap A^c\cap Y^c=\]
\[=\big(M_C\cap Q_C\cap A^c\cap Y^c\big)\cup\big(Z\cap Q_C\cap A^c\cap Y^c\big)\sub\]
\[\sub(M_C\cap A^c)\cup(Z\cap Y^c).\]
Again, by (C.3), we have
\[\tag{D.4}\lambda(Z\cap Y^c)\le\eta/30n.\]
Since $M_C\cap P=\emptyset$, by the definition of $A$, for $M_C\cap A^c$ we have
\[M_C\cap A^c=M_C\cap(X\cup H)\sub\tau(X\cup H)^c\cap(X\cup H)=(X\cup H)\sm\tau(X\cup H),\]
hence, by (C.3), we get
\[\tag{D.5}\lambda(M_C\cap A^c)\le\eta/30n.\]
Inequalities (D.3)--(D.5) yield
\[\tag{E.2}\lambda\big(((M_C\cup Z)\cap Q_C)\sm((A\cup Y)\cap C)\big)\le3\eta/30n.\]

By (C.5), (E.1), and (E.2), we have
\[\tag{C.6}\lambda\big(((A\cup Y)\cap C)\triangle(M_C\cup(Z\cap Q_C))\big)\le5\eta/30n.\]

%
%
Recall that $R=\tau(W)$, so $R\cap Q_C=\tau(W)\cap\tau(C)=\tau(W\cap C)$, hence, by (C.3), we have
\[\tag{C.7}\lambda((W\cap C)\triangle(R\cap Q_C))\le\eta/30n.\]

By Lemma \ref{lemma:ineqpsi2} and inequalities (C.4), (C.6), and (C.7), it holds
\[\psi_m((A\cup Y)\cap C),\ W\cap C)\le\psi_m(M_C\cup(Z\cap Q_C),\ R\cap Q_C)+\]
\[+\lambda\big(((A\cup Y)\cap C)\triangle(M_C\cup(Z\cap Q_C))\big)+\lambda((W\cap C)\triangle(R\cap Q_C))\le\]
\[\le\eta/(|\dD_0|\cdot m)+5\eta/30n+\eta/30n\le\eta/2m+\eta/5n\le\eta/2m+\eta/2m=\eta/m,\]
which proves that condition (A.3) holds for $t<m\le n$, too.
%
%
%
%
\end{proof}

Before we proceed with the next lemma we need to introduce an auxiliary notion. Recall that $\theta$ is the pointwise limit of the sequence $\seqn{|\nu_n|(\cdot)}$.

\begin{definition}\label{def:good}
Fix $p'\io$. Let $\seq{\bB_q}{0\le q\le p'}$ be an increasing sequence of subsets of $\bB$, $\seq{N_q}{0\le q\le p'}$ a sequence of pairwise disjoint elements of $\cC$, and $\seq{m_q}{0\le q\le p'}$ and $\seq{n_q}{0\le q\le2p'+1}$ strictly increasing sequences of natural numbers. For every $0\le p\le p'$ put
\[B_p=\bigcup_{q=0}^pN_q.\]
Let us say that the quadruple
\[\Big\langle\seq{\bB_q}{q\le p'},\ \seq{N_q}{q\le p'},\ \seq{m_q}{q\le p'},\ \seq{n_q}{q\le2p'+1}\Big\rangle\]
is \textbf{\textit{good}} if for every $0\le p\le p'$ the following conditions are satisfied:
\begin{enumerate}[(G.1)]
	\item $N_p\cap H_{n_r}=\emptyset$ for every $r<2p+1$,
	\item $|\nu_{n_{2p+1}}(N_p\cap H_{n_{2p+1}})|\ge0.4$,
	\item 
	\begin{enumerate}[(a)]
		\item $\theta(B_p)<0.1$,
	\end{enumerate}
	and if $p>0$, then:
	\begin{enumerate}[(a),resume]	
		\item $|\nu_{n_{2p}}|(B_{p-1})<0.1$,
		\item $|\nu_{n_{2p+1}}|(B_{p-1})<0.1$,
	\end{enumerate}
	\item for every $0\le q\le p$, $A\in\bB_q$, and atom $U$ of $\aA_{m_q}$, the following hold:
	\begin{enumerate}[(a)]
		\item either
		\[\lambda((A\cap B_p)\cap U)\le2^{-q}\lambda(U)(1-2^{-p-1})/m_q,\]
		or
		\[\lambda(U\sm (A\cap B_p))\le2^{-q}\lambda(U)(1-2^{-p-1})/m_q,\]
		\item for every $m>m_q$ we have
		\[\varphi_m((A\cap B_p)\cap U)\le2^{-q}\lambda(U)(1-2^{-p-1})/m,\]
		\item either
		\[\lambda(A\cap U)\le2^{-q}\lambda(U)/m_q,\]
		or
		\[\lambda(U\sm A)\le2^{-q}\lambda(U)/m_q,\]
		\item for every $m>m_q$ we have
		\[\varphi_m(A\cap U)\le2^{-q}\lambda(U)/m.\]
	\end{enumerate}
\end{enumerate}
\end{definition}

\medskip

\begin{remark}
Trivially, if a quadruple
\[\Big\langle\seq{\bB_q}{q\le p'},\ \seq{N_q}{q\le p'},\ \seq{m_q}{q\le p'},\ \seq{n_q}{q\le2p'+1}\Big\rangle\]
is good, then for any $0\le p\le p'$ the quadruple
\[\Big\langle\seq{\bB_q}{q\le p},\ \seq{N_q}{q\le p},\ \seq{m_q}{q\le p},\ \seq{n_q}{q\le2p+1}\Big\rangle\]
is also good.
\end{remark}

\begin{lemma}
\label{lemma:deuxieme}
Fix $p\io$. Let a quadruple
\[e=\Big\langle\seq{\bB_q}{q\le p},\ \seq{N_q}{q\le p},\ \seq{m_q}{q\le p},\ \seq{n_q}{q\le2p+1}\Big\rangle,\]
be good. Assume that the set $\bB_p$ is finite and let $B\in\bB\sm\bB_p$. Set $\bB_{p+1}=\bB_p\cup\{B\}$. Then, there is a set $N_{p+1}\in\cC$ disjoint with $\bigcup_{q=0}^pN_q$ and natural numbers $m_{p+1}>m_p$ and $n_{2p+3}>n_{2p+2}>n_{2p+1}$ such that the quadruple
\[e'=\Big\langle\seq{\bB_q}{q\le p+1},\ \seq{N_q}{q\le p+1},\ \seq{m_q}{q\le p+1},\ \seq{n_q}{q\le2p+3}\Big\rangle\]
is good.
\end{lemma}
\begin{proof}
%
%
%
Recall that
\[B_p=\bigcup_{q=0}^pN_q.\]
Since $\bB_{p+1}$ is finite and $\bB$ has property (T), by condition (T.1) (applied to $\eps=2^{-p-2}$), there is $m_{p+1}>m_p$ such that for every set $A\in\bB_{p+1}$, atom $U$ of $\aA_{m_{p+1}}$, and set $A'\in\{A,A\cap B_p\}$, the following hold:
\begin{enumerate}[{(F.1)}]
	\item either
	\[\lambda(A'\cap U)\le2^{-p-2}\lambda(U)/m_{p+1},\]
	or
	\[\lambda(U\sm A')\le2^{-p-2}\lambda(U)/m_{p+1},\]
	\item for every $m>m_{p+1}$ we have
	\[\varphi_m(A'\cap U)\le2^{-p-2}\lambda(U)/m.\]
\end{enumerate}

\medskip

Let $\dD_0$ denote the set of all atoms of the finite subalgebra $\dD$ of $\bB$ generated by
\[\bB_{p+1}\cup\{B_p\}\cup\{H_0,\ldots,H_{n_{2p+1}}\}\cup\aA_{m_{p+1}},\]
and define
\[\dD_0^+=\{C\in\dD_0\colon\ \lambda(C)>0\}\]
and
\[\eps=\frac{1}{100}\min\{\lambda(C)\colon\ C\in\dD_0^+\}.\]
Note that $m_{p+1}\ge1$, so $|\dD_0^+|\ge2$. By condition (T.1) there exists $t\io$, $t>m_{p+1}$, such that for every $C\in\dD$ and atom $U$ of $\aA_t$ we have either $\lambda(C\cap U)<\eps\lambda(U)/t$, or $\lambda(U\sm C)<\eps\lambda(U)/t$. Consequently, for every $C\in\dD_0^+$ there is an atom $U_C$ of $\aA_t$ such that
\[\lambda(U_C\cap C)>0.99\lambda(U_C).\]
To see this, note that there exists an atom $U'$ of $\aA_t$ such that $\lambda(U'\sm C)<\eps\lambda(U')/t$, otherwise, by the definition of $\eps$, we would have
\[\lambda(C)=\sum_{U\text{ -- atom of }\aA_t}\lambda(C\cap U)\le\eps/t\cdot\sum_{U\text{ -- atom of }\aA_t}\lambda(U)=\eps/t<\lambda(C),\]
which is impossible. Since $\eps\le 1/100$, for such $U'$ we have
\[\lambda(U'\cap C)=\lambda(U')-\lambda(U'\sm C)>\lambda(U')-\eps\lambda(U')/t\ge\lambda(U')-\lambda(U')/100=0.99\lambda(U').\]
Set $U_C=U'$.

\medskip

Let $\xi=0.1-\theta(B_p)$. By condition (G.3.a), $\xi>0$. Let also $\delta>0$ be such that for every $A\in\cC$ if $\lambda(A)<\delta$, then $\theta_1(A)\le\xi/3$ (the existence of such $\delta$ follows from the Lebesgue decomposition of $\theta$, see the beginning of this section).

Set
\[\gamma=\min\big(\delta/2,\ 2^{-m_{p+1}-2p-5}/t,\ 2^{-t-12}\big)\]
and
\[\zeta=2^{-10}\gamma^2/|\dD_0|^4.\]
Since $\gamma<1$, it holds $\zeta<\gamma/3\le\delta/2$. Let $n_{2p+3}>n_{2p+2}>n_{2p+1}$ be such integers that
\[\lambda(H_{n_{2p+2}})<\zeta/2,\quad\lambda(H_{n_{2p+3}})<\zeta/2,\quad\theta_2(H_{n_{2p+3}})<\xi/3,\]
as well as
\[|\nu_{n_{2p+2}}|(B_p)<0.1\quad\text{and}\quad|\nu_{n_{2p+3}}|(B_p)<0.1\]
(for the last two inequalities recall that the sequence $\seqn{|\nu_n|(\cdot)}$ converges to $\theta$ pointwise). Since $|\nu_{n_{2p+3}}|(H_{n_{2p+3}})\ge0.95$, we have $|\nu_{n_{2p+3}}|(H_{n_{2p+3}}\sm B_p)\ge0.85$. Consequently, 
there exists $Y\in\cC$ such that $Y\sub H_{n_{2p+3}}\sm B_p$ and $|\nu_{n_{2p+3}}(Y)|\ge0.4$. We also have
\[\tag{F.3}\lambda(Y)\le\lambda(H_{n_{2p+3}})<\zeta<\gamma\le\delta/2.\]
Let $X\in\cC$ be such that $\lambda(X)<\zeta$ and $\theta_2(X^c)<\xi/3$ (the existence of such $X$ follows again from the Lebesgue decomposition of $\theta$). 

\medskip

Let a set $A\in\cC$ be as in Lemma \ref{lemma:troisieme} (applied to $\dD_0=\dD_0^+$, $t$, $\eta=\gamma/|\dD_0|$, $\zeta$, $X$, $W=B_p\cap(H_{n_{2p+2}}\cup H_{n_{2p+3}})$, $Y$, $H=H_{n_{2p+2}}\cup H_{n_{2p+3}}$). Put
\[N=A\sm(B_p\cup H_0\cup\ldots\cup H_{n_{2p+1}}).\]
Then, $N\in\cC$. By condition (A.1) we have
\[\tag{F.4}N\cap(X\cup H_{n_{2p+2}}\cup H_{n_{2p+3}})=\emptyset,\]
and, by condition (A.2),
\[\tag{F.5}\lambda(N)\le\gamma/|\dD_0|\le\gamma\le\delta/2.\]

Let $C\in\dD_0$ and $m>t$. We will estimate $\varphi_m((N\cup Y)\cap C)$. If $\lambda(C)=0$, then $\varphi_m((N\cup Y)\cap C)=0$. Assume thus that $C\in\dD_0^+$. If $C\sub B_p\cup(H_0\cup\ldots\cup H_{n_{2p+1}})$, then by the definition of $N$ and the fact that $Y\sub H_{n_{2p+3}}\sm B_p$, we get that $(N\cup Y)\cap C=\emptyset$, so
\[\varphi_m((N\cup Y)\cap C)=\varphi_m(\emptyset)=0\le\gamma/(m\cdot|\dD_0|).\]
If, on the other hand, $C\not\subseteq B_p\cup H_0\cup\ldots\cup H_{n_{2p+1}}$, then $C\sub(B_p\cup H_0\cup\ldots\cup H_{n_{2p+1}})^c$, as $C$ is an atom of $\dD$. We thus have
\[(A\cup Y)\cap C=(A\cap C)\cup(Y\cap C)=\big(A\cap(B_p\cup H_0\cup\ldots\cup H_{n_{2p+1}})^c\cap C\big)\cup(Y\cap C)=\]
\[=(N\cap C)\cup(Y\cap C)=(N\cup Y)\cap C.\]
Hence, since $B_p\cap C=\emptyset$ and by condition (A.3), we get
\[\varphi_m((N\cup Y)\cap C)=\varphi_m((A\cup Y)\cap C)=\psi_m((A\cup Y)\cap C), \emptyset)=\]
\[=\psi_m\big((A\cup Y)\cap C,\ B_p\cap (H_{n_{2p+2}}\cup H_{n_{2p+3}})\cap C\big)\le\gamma/(m\cdot|\dD_0|).\]
In both cases we get that for every $C\in\dD_0$ and $m>t$ it holds
\[\tag{F.6}\varphi_m((N\cup Y)\cap C)\le\gamma/(m\cdot|\dD_0|).\]

\medskip

Put $N_{p+1}=N\cup Y$. Of course, $N_{p+1}\in\cC$. Since $N\cap B_p=\emptyset$ and $Y\sub H_{n_{2p+3}}\sm B_p$, $N_{p+1}\cap B_p=\emptyset$ and hence $N_{p+1}\cap N_q=\emptyset$ for any $q<p+1$. 

By $\lambda(N)\le\delta/2$ and (F.3), we have $\lambda(N_{p+1})<\delta$, so by the choice of $\delta$ it holds
\[\theta_1(N_{p+1})\le\xi/3.\]
Since $N\sub X^c$ and $Y\sub H_{n_{2p+3}}$, we have as well
\[\theta_2(N_{p+1})=\theta_2(N)+\theta_2(Y)<\xi/3+\xi/3=2\xi/3.\]
It follows that $\theta(N_{p+1})<\xi$.

Put
\[B_{p+1}=\bigcup_{q=0}^{p+1}N_q.\]
Because $B_{p+1}=B_p\cup N_{p+1}$, we get that
\[\theta(B_{p+1})=\theta(B_p\cup N_{p+1})=\theta(B_p)+\theta(N_{p+1})<\theta(B_p)+\xi=\theta(B_p)+0.1-\theta(B_p)=0.1,\]
so $\theta(B_{p+1})<0.1$, hence condition (G.3.a) holds for the set $B_{p+1}$. Conditions (G.3.b) and (G.3.c) hold by the choice of $n_{2p+2}$ and $n_{2p+3}$, so, consequently, condition (G.3) is satisfied by the quadruple 
\[e'=\Big\langle\seq{\bB_q}{q\le p+1},\ \seq{N_q}{q\le p+1},\ \seq{m_q}{q\le p+1},\ \seq{n_q}{q\le2p+3}\Big\rangle.\]
Since $Y\sub H_{n_{2p+3}}$, the definition of $N$ and equality (F.4) imply that $N_{p+1}\cap H_{n_r}=\emptyset$ for every $r<2p+3$, that is, condition (G.1) holds for the quadruple $e'$. Condition (G.2) is also satisfied by $e'$, because $N_{p+1}\cap H_{n_{2p+3}}=Y$ and $|\nu_{n_{2p+3}}(Y)|\ge0.4$. 


\medskip

To finish the proof, we show that condition (G.4) holds for $e'$, too. Fix $q\le p+1$, $A\in\bB_q$, and an atom $U$ of $\aA_{m_q}$. If $q\le p$, then conditions (G.4.c) and (G.4.d) follow from the fact that the quadruple $e$ is good, and if $q=p+1$, then these conditions follow directly from conditions (F.1) and (F.2) (with $A'=A$). Consequently, $e'$ satisfies conditions (G.4.c) and (G.4.d). We thus need only to prove that $e'$ satisfies conditions (G.4.a) and (G.4.b).

We first show that $e'$ satisfies condition (G.4.a). As $B_{p+1}\sm B_p=N_{p+1}$ and $B_p\sm B_{p+1}=\emptyset$, we have $B_{p+1}\triangle B_p=N_{p+1}$ and so, by inequalities (F.3) and (F.5), we get
\[\lambda(B_{p+1}\triangle B_p)=\lambda(N_{p+1})=\lambda(N\cup Y)=\lambda(N)+\lambda(Y)<2\gamma,\]
By the definition of $\gamma$, we have
\[\tag{F.7}\lambda(B_{p+1}\triangle B_p)\le2\cdot2^{-m_{p+1}-2p-5}/t=2^{-m_{p+1}-2p-4}/t=2^{-(p+1)}2^{-m_{p+1}-1}2^{-(p+2)}/t.\]

Assume that $q\le p$. By condition (G.4.a) (for $e$) either
\[\lambda((A\cap B_p)\cap U)\le2^{-q}\lambda(U)(1-2^{-p-1})/m_q,\]
or
\[\lambda(U\sm(A\cap B_p))\le2^{-q}\lambda(U)(1-2^{-p-1})/m_q.\]
If
\[\lambda((A\cap B_p)\cap U)\le2^{-q}\lambda(U)(1-2^{-p-1})/m_q,\]
then, by Lemma \ref{lemma:ineq}.(1) (with $B=B_p$ and $B'=B_{p+1}$) and (F.7), we have
\[\lambda((A\cap B_{p+1})\cap U)\le\lambda((A\cap B_p)\cap U)+\lambda(B_p\triangle B_{p+1})\le\]
\[\le2^{-q}\lambda(U)(1-2^{-p-1})/m_q+2^{-(p+1)}2^{-m_{p+1}-1}2^{-(p+2)}/t\le\]
\[\le2^{-q}\lambda(U)(1-2^{-p-1})/m_q+2^{-q}\lambda(U)2^{-(p+2)}/m_q=2^{-q}\lambda(U)(1-2^{-(p+2)})/m_q.\]
If, on the other hand, it holds
\[\lambda(U\sm(A\cap B_p))\le2^{-q}\lambda(U)(1-2^{-p-1})/m_q,\]
then, by Lemma \ref{lemma:ineq}.(2) (with $B=B_p$ and $B'=B_{p+1}$) and (F.7), we similarly have
\[\lambda(U\sm(A\cap B_{p+1}))\le\lambda(U\sm(A\cap B_p))+\lambda(B_{p+1}\triangle B_p)\le\]
\[\le2^{-q}\lambda(U)(1-2^{-p-1})/m_q+2^{-q}\lambda(U)2^{-(p+2)}/m_q=2^{-q}\lambda(U)(1-2^{-(p+2)})/m_q.\]

Assume now that $q=p+1$. By Lemma \ref{lemma:ineq}.(1) (with $B=B_p$ and $B'=B_{p+1}$) and conditions (F.1) (with $A'=A\cap B_p$) and (F.7), we either have
\[\lambda((A\cap B_{p+1})\cap U)\le\lambda((A\cap B_p)\cap U)+\lambda(B_p\triangle B_{p+1})\le\]
\[\le2^{-p-2}\lambda(U)/m_{p+1}+2^{-(p+1)}2^{-m_{p+1}-1}2^{-(p+2)}/t\le2^{-(p+1)}\lambda(U)(2^{-1}+2^{-(p+2)})/m_{p+1}\le\]
\[\le2^{-(p+1)}\lambda(U)(1-2^{-(p+2)})/m_{p+1},\]
or, by similar calculations, it holds
\[\lambda(U\sm(A\cap B_{p+1}))\le\lambda(U\sm(A\cap B_p))+\lambda(B_{p+1}\triangle B_p)\le\]
\[\le2^{-p-2}\lambda(U)/m_{p+1}+2^{-(p+1)}2^{-m_{p+1}-1}2^{-(p+2)}/t\le2^{-(p+1)}\lambda(U)(1-2^{-(p+2)})/m_{p+1}.\]
It follows that $e'$ satisfies condition (G.4.a).


\medskip

Finally, we prove that condition (G.4.b) holds for $e'$, too---this will finish the proof of the lemma. 
Fix $m>m_q$. We have here four cases in total.

If $q\le p$ and $m\le t$, then, by Lemma \ref{lemma:ineq}.(5) (with $B=B_p$ and $B'=B_{p+1}$) and conditions (G.4.b) (for $e$) and (F.7), we have
\[\varphi_m((A\cap B_{p+1})\cap U)\le\varphi_m((A\cap B_p)\cap U)+\lambda(B_{p+1}\triangle B_p)\le\]
\[\le2^{-q}\lambda(U)(1-2^{-p-1})/m+2^{-(p+1)}2^{-m_{p+1}-1}2^{-(p+2)}/t\le\]
\[\le2^{-q}\lambda(U)(1-2^{-p-1})/m+2^{-q}\lambda(U)2^{-(p+2)}/m=2^{-q}\lambda(U)(1-2^{-(p+2)})/m.\]

If $q=p+1$ and $m\le t$, then, by Lemma \ref{lemma:ineq}.(5) (with $B=B_p$ and $B'=B_{p+1}$) and conditions (F.2) (with $A'=A\cap B_p$) and (F.7), we have
\[\varphi_m((A\cap B_{p+1})\cap U)\le\varphi_m((A\cap B_p)\cap U)+\lambda(B_{p+1}\triangle B_p)\le\]
\[\le2^{-p-2}\lambda(U)/m+2^{-(p+1)}2^{-m_{p+1}-1}2^{-(p+2)}/t\le2^{-(p+1)}\lambda(U)(2^{-1}+2^{-(p+2)})/m\le\]
\[\le2^{-(p+1)}\lambda(U)(1-2^{-(p+2)})/m.\]

Assume now that $m>t$---we deal with this case a bit differently. 
Since $A\cap U\in\dD$ is the union of at most $|\dD_0|$ many atoms of the algebra $\dD$, by Lemma \ref{lemma:ineqpsi}.(1) (with $W=N_{p+1}$ and sets $Z$ and $Z'$ ranging over all atoms $C\in\dD_0$), inequality (F.6), and the definition of $\gamma$, we have
\[\varphi_m(N_{p+1}\cap(A\cap U))\le|\dD_0|\cdot\gamma/(m\cdot|\dD_0|)\le\]
\[\le2^{-m_{p+1}-2p-5}/(t\cdot m)\le2^{-(p+1)}\lambda(U)2^{-(p+2)}/m.\]
It follows from Lemma \ref{lemma:ineq}.(3) (with $B=B_p$ and $B'=B_{p+1}$) and the above inequality that
\[\tag{F.8}\varphi_m((A\cap B_{p+1})\cap U)\le\]
\[\le\varphi_m((A\cap B_p)\cap U)+\psi_m\big((A\cap(B_{p+1}\sm B_p)\cap U,\ (A\cap (B_p\sm B_{p+1}))\cap U\big)=\]
\[=\varphi_m((A\cap B_p)\cap U)+\psi_m(N_{p+1}\cap(A\cap U),\emptyset)=\varphi_m((A\cap B_p)\cap U)+\varphi_m(N_{p+1}\cap(A\cap U))\le\]
\[\le\varphi_m((A\cap B_p)\cap U)+2^{-(p+1)}\lambda(U)2^{-(p+2)}/m.\]

If $q\le p$, then, by (F.8) and condition (G.4.b) (for $e$), we have
\[\varphi_m((A\cap B_{p+1})\cap U)\le\varphi_m((A\cap B_p)\cap U)+2^{-(p+1)}\lambda(U)2^{-(p+2)}/m\le\]
\[\le2^{-q}\lambda(U)(1-2^{-p-1})/m+2^{-(p+1)}\lambda(U)2^{-(p+2)}/m\le\]
\[\le2^{-q}\lambda(U)(1-2^{-p-1})/m+2^{-q}\lambda(U)2^{-(p+2)}/m=2^{-q}\lambda(U)(1-2^{-(p+2)})/m.\]

If $q=p+1$, then, by conditions (F.8) and (F.2) (with $A'=A\cap B_p$), we have
\[\varphi_m((A\cap B_{p+1})\cap U)\le\varphi_m((A\cap B_p)\cap U)+2^{-(p+1)}\lambda(U)2^{-(p+2)}/m\le\]
\[\le2^{-(p+2)}\lambda(U)/m+2^{-(p+1)}\lambda(U)2^{-(p+2)}/m=2^{-(p+1)}\lambda(U)(2^{-1}+2^{-(p+2)})/m\le\]
\[\le2^{-(p+1)}\lambda(U)(1-2^{-(p+2)})/m.\]
It follows that condition (G.4.b) holds also for $e'$. Consequently, the quadruple $e'$ satisfies condition (G.4) and ultimately it is a good quadruple. The proof of the lemma is thus finished.
\end{proof}

\begin{remark}\label{remark:condition_a3}
Notice that condition (A.3) from the conclusion of Lemma \ref{lemma:troisieme}, satisfied by the set $A$, was relevant in the proof of Lemma \ref{lemma:deuxieme} only in the case of atoms $C\in\dD_0^+$ such that $B_p\cap(H_{n_{2p+2}}\cup H_{n_{2p+3}})\cap C=\emptyset$.
\end{remark}

\begin{remark}\label{remark:good_conditions_exist}
Let us too note that the argument presented in the proof of Lemma \ref{lemma:deuxieme} works also if we start with the ``empty'' quadruple $e=\big\langle\langle\rangle,\langle\rangle,\langle\rangle,\langle\rangle\big\rangle$ and any set $B\in\bB$. This proves that there exist good quadruples of any positive length (cf. the density of the sets $E_n$ in the next proof).
\end{remark}

The proof of the next lemma is actually the only place in this paper where Martin's axiom for $\sigma$-centered posets, $\mathsf{MA}(\sigma\text{-centered})$, is invoked. Recall that a poset $(\PP,\le)$ is \textit{$\sigma$-centered} if $\PP$ can be written as a countable union $\PP=\bigcup_{i\io}\PP_i$, where each set $\PP_i$ is \textit{centered}, that is, for every $e_1,\ldots,e_n\in\PP_i$, $n\io$, $n>1$, there is $e\in\PP$ such that $e\le e_j$ for each $j=1,\ldots,n$.

\begin{lemma}
\label{lemma:intermediaire}
Assume that $\mathsf{MA}(\sigma\text{-centered})$ holds. If $|\bB|<\frakc$, then there exist an increasing sequence $\seqp{\bB_p}$ of subsets of $\bB$ such that $\bB=\bigcup_{p\io}\bB_p$, a sequence $\seqp{N_p}$ of pairwise disjoint elements of $\cC$, and strictly increasing sequences $\seqp{m_p\io}$ and $\seqp{n_p\io}$ such that for every $p'\io$ the quadruple \[\Big\langle\seq{\bB_p}{p\le p'},\ \seq{N_p}{p\le p'},\ \seq{m_p}{p\le p'},\  \seq{n_p}{p\le2p'+1}\Big\rangle\]
is good.
\end{lemma}
\begin{proof}
Assume that $|\bB|<\frakc$. We define a $\sigma$-centered poset $(\PP,\le)$ as follows. Conditions of $\PP$ are all good quadruples
\[e=\big\langle\ol{\bB}_e,\ \ol{N}_e,\ \ol{m}_e,\ \ol{n}_e\big\rangle,\]
consisting of finite sequences such that, for some number $p^e\io$ depending on $e$, we have:
\begin{itemize}
	\item $\ol{\bB}_e=\seq{\bB_p^e}{p\le p^e}$ is an increasing sequence of finite\footnote{Note, however, that in the definition of a good pair (Definition \ref{def:good}) we do note require that the sets $\bB_q$ are finite.} subsets of $\bB$,
	\item $\ol{N}_e=\seq{N_p^e}{p\le p^e}$ is a sequence of pairwise disjoint elements of $\cC$,
	\item $\ol{m}_e=\seq{m_p^e}{p\le p^e}$ and $\ol{n}_e=\seq{n_p^e}{p\le 2p^e+1}$ are strictly increasing sequences of natural numbers.
\end{itemize}
By Remark \ref{remark:good_conditions_exist}, $\PP\neq\emptyset$. For $e_1,e_2\in\PP$, where $e_1=\big\langle\ol{\bB}_{e_1},\ol{N}_{e_1},\ol{m}_{e_1},\ol{n}_{e_1}\big\rangle$ and $e_2=\big\langle\ol{\bB}_{e_2},\ol{N}_{e_2},\ol{m}_{e_2},\ol{n}_{e_2}\big\rangle$, we write
\[e_1\ge e_2\]
and say that $e_2$ \textit{extends} $e_1$ if $p^{e_1}\le p^{e_2}$ and for every $q\le p^{e_1}$ we have: $\bB_q^{e_1}\sub\bB_q^{e_2}$, $N_q^{e_1}=N_q^{e_2}$, $m_q^{e_1}=m_q^{e_2}$, and $n_{2q}^{e_1}=n_{2q}^{e_2}$ and $n_{2q+1}^{e_1}=n_{2q+1}^{e_2}$.

\medskip

We claim that the poset $(\PP,\le)$ is $\sigma$-centered. Let
\[\varphi\colon\N\to\N\times\cC^{<\omega}\times\N^{<\omega}\times\N^{<\omega}\]
be a bijection (recall that $\cC$ is countable)\footnote{Here, $A^{<\omega}$ denotes the family of all finite sequences consisting of elements of a set $A$.}. For each $i\io$ let $\PP_i$ consist of all those conditions $e=\big\langle \ol{\bB}_e,\ol{N}_e,\ol{m}_e,\ol{n}_e\big\rangle\in\PP$ such that
\[\varphi(i)=\big\langle p^e,\ol{N}_e,\ol{m}_e,\ol{n}_e\big\rangle.\]
Then, $\PP=\bigcup_{i\io}\PP_i$, as $\varphi$ is a bijection. Fix $i\io$ and let $e_1,\ldots,e_n\in\PP_i$ for some $n\io$, $n>1$. Set
\[e=\Big\langle\seq{\bB_q^{e_1}\cup\ldots\cup\bB_q^{e_n}}{q\le p^{e_1}},\ \ol{N}_{e_1},\ \ol{m}_{e_1},\ \ol{n}_{e_1}\Big\rangle.\]
We first need to show that $e$ is good. We have $p^e=p^{e_1}$. Let $p\le p^e$. Conditions (G.1)--(G.3) are independent of the first coordinate for any good pair, hence they are trivially satisfied for $p$ and $e$ as $e_1$ is good. Let $q\le p$, $A\in\bB_q^{e_1}\cup\ldots\cup\bB_q^{e_n}$, and $U$ be an atom of $\aA_{m_q}$. There is $1\le j\le n$ such that $A\in\bB_q^{e_j}$. Since $e_j$ is good, conditions (G.4.a)--(G.4.d) are satisfied for $q$, $p$, $A$, and $U$. Consequently, $e$ satisfies condition (G.4) and hence $e$ is good, too. It follows that $e\in\PP_i$. Since $e$ extends each of the conditions $e_1,\ldots,e_n$, we get that $\PP_i$ is centered. Consequently, $\PP$ is $\sigma$-centered. 

\medskip

For every $B\in\bB$ put
\[D_B=\{e\in\PP\colon\ B\in\bB_{p^e}^e\}.\]
Fix $B\in\bB$. For every $e\in\PP\sm D_B$ by Lemma \ref{lemma:deuxieme} there is $e'\in\PP$ extending $e$ and such that $p^{e'}=p^e+1$ and $\bB_{p^{e'}}^{e'}=\bB_{p^e}^e\cup\{B\}$, which implies that $e'\in D_B$. This shows that the set $D_B$ is dense in $\PP$. Similarly, for every $n\io$, the set
\[E_n=\{e\in\PP\colon\ p^e\ge n\}\]
is dense in $\PP$. Since $|\bB|<\frakc$, the collection $\{D_B\colon B\in\bB\}\cup\{E_n\colon n\io\}$ has cardinality $<\frakc$, thus, by $\mathsf{MA}(\sigma\text{-centered})$, there exists a filter $G\sub\PP$ such that $G\cap D_B\neq\emptyset$ for every $B\in\bB$ and $G\cap E_n\neq\emptyset$ for every $n\io$. It follows that
\[\bB=\bigcup_{e\in G}\bB_{p^e}^e.\]
For each $p\io$ put
\[\bB_p=\bigcup\{\bB_p^e\colon\ e\in G,\ p\le p^e\}.\]

The sequence $\seqp{\bB_p}$ is increasing. Indeed, let $p\io$ and $B\in\bB_p$. There is $e\in G$ such that $p\le p^e$ and $B\in\bB_p^e$. By the density of $E_{p+1}$, there is also $e'\in G$ such that $p+1\le p^{e'}$. Since $G$ is a filter, there is $e''\in G$ such that $e''\le e$ and $e''\le e'$. Consequently, $B\in\bB_p^e\sub\bB_p^{e''}$ and $p+1\le p^{e''}$, so $B\in\bB_{p+1}^{e''}$, and hence $B\in\bB_{p+1}$. It follows that $\bB_p\sub\bB_{p+1}$.

For every $B\in\bB$ there is $e\in G\cap D_B$, so $B\in\bB_{p^e}^e$, hence $B\in\bB_{p^e}$ and ultimately $B\in\bigcup_{p\io}\bB_p$. We thus have $\bB=\bigcup_{p\io}\bB_p$.

Set also:
\[\seqp{N_p}=\bigcup_{\substack{e\in G}}\ol{N}_e,\quad\seqp{m_p}=\bigcup_{\substack{e\in G}}\ol{m}_e,\quad\seqp{n_p}=\bigcup_{\substack{e\in G}}\ol{n}_e.\]
Of course, since $G$ is a filter, the sets $N_p$'s are pairwise disjoint and the sequences $\seqp{m_p}$ and $\seqp{n_p}$ are strictly increasing.

Fix $p'\io$. We need to show that the quadruple
\[e=\Big\langle\seq{\bB_p}{p\le p'},\ \seq{N_p}{p\le p'},\ \seq{m_p}{p\le p'},\  \seq{n_p}{p\le2p'+1}\Big\rangle\]
is good. Let $p\le p'$ and $q\le p$. Let $A\in\bB_q$ and $U$ be an atom of $\aA_{m_q}$. There is $e'\in G$ such that $q\le p^{e'}$ and $A\in\bB_q^{e'}$. Since $G$ is a filter, there is also $e''\in G$ such that $p^{e''}\ge p'$ and for every $r\le p'$ we have: $N_r=N_r^{e''}$, $m_r=m_r^{e''}$, and $n_{2r}=n_{2r}^{e''}$ and $n_{2r+1}=n_{2r+1}^{e''}$. Let $f\in G$ be such that $f\le e'$ and $f\le e''$. It follows that $A\in\bB_q^f$, $p\le p'\le p^f$, and for every $r\le p$ we have: $N_r=N_r^f$, $m_r=m_r^f$, and $n_{2r}=n_{2r}^f$ and $n_{2r+1}=n_{2r+1}^f$. Since $f$ is good and $p\le p^f$, conditions (G.1)--(G.3) are satisfied for $p$ and $f$, and so for $p$ and $e$. Similarly, since $f$ is good and $q\le p\le p^f$ and $A\in\bB_q^f$, conditions (G.4.a)--(G.4.d) hold for $q$, $p$, $A$, and $U$, and hence $e$ satisfies condition (G.4) as well. Consequently, $e$ satisfies all of conditions (G.1)--(G.4) for any $p\le p'=p^e$ and hence $e$ is good, too. The proof of the lemma is thus finished.
\end{proof}

\begin{lemma}
\label{lemma:premiere}
%
Assume that $\mathsf{MA}(\sigma\text{-centered})$ holds. If $|\bB|<\frakc$, then there exist a set $B\in Bor(\Cantor)$ and a strictly increasing sequence $\seqk{n_k\io}$ such that the Boolean algebra $\tilde{\bB}$ generated by $\bB\cup\{B\}$ has property (T) and for every $k\io$ it holds:
\begin{enumerate}[({B}.1)]
	\item $B\cap H_{n_k}\in\cC$,
	\item $|\nu_{2k}(B\cap H_{n_{2k}})|\le0.1$,
	\item $|\nu_{2k+1}(B\cap H_{n_{2k+1}})|\ge0.3$.
\end{enumerate}
\end{lemma}
\begin{proof}

Assume that $|\bB|<\frakc$. Let $\seqp{\bB_p}$, $\seqp{N_p}$, $\seqp{m_p}$, $\seqp{n_p}$, and $\seqp{B_p}$ be as in Lemma \ref{lemma:intermediaire}. Put
\[B=\bigcup_{p=0}^\infty B_p=\bigcup_{q=0}^\infty N_q,\]
and let $\tilde{\bB}$ be the Boolean algebra generated by $\bB\cup\{B\}$. We first show that the set $B$ satisfies conditions (B.1)--(B.3). By condition (G.1) for every $k\io$ we have
\[B\cap H_{n_k}=\bigcup_{q=0}^\infty N_q\cap H_{n_k}=\Big(\bigcup_{q\le\frac{k-1}{2}}N_q\Big)\cap H_{n_k},\]
which clearly belongs to $\cC$, as the latter union is finite for fixed $k$ and each $N_q$ as well as $H_{n_k}$ is in $\cC$. Consequently, (B.1) is satisfied.

If $k=0$, then $B\cap H_{n_k}=\emptyset$, hence
\[|\nu_{n_{2k}}(B\cap H_{n_{2k}})|=|\nu_{n_{2k}}(\emptyset)|=0<0.1.\]
If $k>0$, then 
\[B\cap H_{n_{2k}}=\Big(\bigcup_{q\le\frac{2k-1}{2}}N_q\Big)\cap H_{n_{2k}}=\Big(\bigcup_{q<k}N_q\Big)\cap H_{n_{2k}}=B_{k-1}\cap H_{n_{2k}},\]
hence, by condition (G.3.b) (with $p'=p=k$), we have
\[|\nu_{n_{2k}}(B\cap H_{n_{2k}})|=|\nu_{n_{2k}}(B_{k-1}\cap H_{n_{2k}})|\le|\nu_{n_{2k}}|(B_{k-1})<0.1,\]
so (B.2) is satisfied for every $k\io$. Similarly, if $k=0$, then $B\cap H_{n_{2k+1}}=N_k\cap H_{n_{2k+1}}$, so, by condition (G.2), we have
\[|\nu_{n_{2k+1}}(B\cap H_{n_{2k+1}})|=|\nu_{n_{2k+1}}(N_k\cap H_{n_{2k+1}})|\ge0.4>0.3.\]
If $k>0$, then
\[B\cap H_{n_{2k+1}}=\Big(\bigcup_{q\le\frac{2k+1-1}{2}}N_q\Big)\cap H_{n_{2k+1}}=\Big(\bigcup_{q\le k}N_q\Big)\cap H_{n_{2k+1}}=\]
\[=\big(\Big(\bigcup_{q<k}N_q\Big)\cap H_{n_{2k+1}}\big)\cup(N_k\cap H_{n_{2k+1}})=(B_{k-1}\cap H_{n_{2k+1}})\cup(N_k\cap H_{n_{2k+1}}),\]
hence, by conditions (G.2) and (G.3.c), we have
\[|\nu_{n_{2k+1}}(B\cap H_{n_{2k+1}})|=\big|\nu_{n_{2k+1}}\big((B_{k-1}\cap H_{n_{2k+1}})\cup(N_k\cap H_{n_{2k+1}})\big)\big|=\]
\[=|\nu_{n_{2k+1}}(B_{k-1}\cap H_{n_{2k+1}})+\nu_{n_{2k+1}}(N_k\cap H_{n_{2k+1}})|\ge\]
\[\ge|\nu_{n_{2k+1}}(N_k\cap H_{n_{2k+1}})|-|\nu_{n_{2k+1}}(B_{k-1}\cap H_{n_{2k+1}})|\ge\]
\[\ge|\nu_{n_{2k+1}}(N_k\cap H_{n_{2k+1}})|-|\nu_{n_{2k+1}}|(B_{k-1})\ge0.4-0.1=0.3,\]
so (B.3) is satisfied for every $k\io$ as well.

\medskip

To finish the proof of the lemma we need to show that the new Boolean algebra $\tilde{\bB}$ has property (T). Since the elements of $\seqp{N_p}$ are pairwise disjoint, for every $p\io$ it holds
\[B\triangle B_p=\bigcup_{q>p}N_q,\]
so, by the $\sigma$-additivity of $\lambda$,
\[\lambda(B\triangle B_p)=\sum_{q>p}\lambda(N_q),\]
hence, by the inequality $\sum_{p=0}^\infty\lambda(N_p)\le 1$,
\[\tag{$\star$}\lim_{p\to\infty}\lambda(B\triangle B_p)=0.\]

Recall that each element of $\tilde{\bB}$ is of the form $(A\cap B)\cup(A'\sm B)$ for some sets $A,A'\in\bB$. So, fix $r\io$ and let $A_1,A_1',\ldots,A_r,A_r'\in\bB$. Let $\eps>0$. Since $\seqq{\bB_q}$ is an increasing sequence with $\bB=\bigcup_{q\io}\bB_q$, there is $q\io$ such that $A_1,A_1',\ldots,A_r,A_r'\in\bB_q$ and $2^{-q}<\eps/4$.

Fix $i=1,\ldots,r$ for a moment and let $X\in\{A_i,A_i'\}$. By condition (G.4.a), for every atom $U$ of $\aA_{m_q}$, there is an infinite set $I_U\sub\N$ with $\min I_U\ge q$ such that for every $p\in I_U$ we have
\[\lambda((X\cap B_p)\cap U)\le2^{-q}\lambda(U)(1-2^{-p-1})/m_q,\]
or there is an infinite set $J_U\sub\N$ with $\min J_U\ge q$ such that for every $p\in J_U$ we have
\[\lambda(U\sm(X\cap B_p))\le2^{-q}\lambda(U)(1-2^{-p-1})/m_q.\]
Fix an atom $U$ of $\aA_{m_q}$ and assume that the first case holds. By Lemma \ref{lemma:ineq}.(1), for every $p\in I_U$, we have
\[\lambda((X\cap B)\cap U)\le\lambda((X\cap B_p)\cap U)+\lambda(B\triangle B_p)\le\]
\[\le2^{-q}\lambda(U)(1-2^{-p-1})/m_q+\lambda(B\triangle B_p),\]
where, by ($\star$), the bottom expression converges to $2^{-q}\lambda(U)/m_q$ as $p\to\infty$, $p\in I_U$. It follows that
\[\lambda((X\cap B)\cap U)\le2^{-q}\lambda(U)/m_q.\]

Assume now that the second case holds. By Lemma \ref{lemma:ineq}.(2), for every $p\in J_U$, we have
\[\lambda(U\sm(X\cap B))\le\lambda(U\sm(X\cap B_p))+\lambda(B\triangle B_p)\le\]
\[\le2^{-q}\lambda(U)(1-2^{-p-1})/m_q+\lambda(B\triangle B_p),\]
where, again by ($\star$), the bottom expression converges to $2^{-q}\lambda(U)/m_q$ as $p\to\infty$, $p\in J_U$. It follows that
\[\lambda(U\sm(X\cap B))\le2^{-q}\lambda(U)/m_q.\]

\medskip

Summing up, for every $i=1,\ldots,r$ and atom $U$ of $\aA_{m_q}$, if $X\in\{A_i,A_i'\}$, then it holds either
\[\lambda((X\cap B)\cap U)\le2^{-q}\lambda(U)/m_q,\]
or
\[\lambda(U\sm(X\cap B))\le2^{-q}\lambda(U)/m_q.\]
Since, by condition (G.4.c), it too holds either
\[\lambda(A_i'\cap U)\le 2^{-q}\lambda(U)/m_q,\]
or
\[\lambda(U\sm A_i')\le2^{-q}\lambda(U)/m_q,\]
Lemma \ref{lemma:complement}.(a) (with $A=A_i'$) implies that it holds either
\[\lambda((A_i'\sm B)\cap U)\le2\cdot2^{-q}\lambda(U)/m_q,\]
or
\[\lambda(U\sm(A_i'\sm B))\le2\cdot2^{-q}\lambda(U)/m_q.\]
From this we get that for every $i=1,\ldots,r$ and atom $U$ of $\aA_{m_q}$ we either have (recall that $2^{-q}<\eps/4$)
\[\lambda\big(((A_i\cap B)\cup(A_i'\sm B))\cap U\big)\le 3\cdot2^{-q}\lambda(U)/m_q<\eps\lambda(U)/m_q,\]
or
\[\lambda\big(U\sm((A_i\cap B)\cup(A_i'\sm B))\big)\le 2\cdot2^{-q}\lambda(U)/m_q<\eps\lambda(U)/m_q.\]
It follows that the algebra $\tilde{\bB}$ satisfies condition (T.1) of property (T).

\medskip

To show that $\tilde{\bB}$ satisfies also condition (T.2) of property (T), fix again $i=1,\ldots,r$, a set $X\in\{A_i,A_i'\}$, and an atom $U$ of $\aA_{m_q}$ and notice that for every $m>m_q$ and $p\ge q$, by Lemma \ref{lemma:ineq}.(5) and condition (G.4.b), we have
\[\varphi_m((X\cap B)\cap U)\le\varphi_m((X\cap B_p)\cap U)+\lambda(B\triangle B_p)\le\]
\[\le2^{-q}\lambda(U)(1-2^{-p-1})/m+\lambda(B\triangle B_p),\]
where, by ($\star$), the bottom expression converges to $2^{-q}\lambda(U)/m$ as $p\to\infty$. Thus, for every $m>m_q$ we have
\[\varphi_m((X\cap B)\cap U)\le2^{-q}\lambda(U)/m.\]
Since, by condition (G.4.d), for every $m>m_q$ we have
\[\varphi_m(A_i'\cap U)\le2^{-q}\lambda(U)/m,\]
by Lemma \ref{lemma:complement}.(b) (with $A=A_i'$), we also get that for every $m>m_q$ it holds
\[\varphi_m((A_i'\sm B)\cap U)\le2\cdot2^{-q}\lambda(U)/m.\]
By Lemma \ref{lemma:ineqpsi}.(1) (with $W=U$, $Z=A_i\cap B$, and $Z'=A_i'\sm B$), the latter two inequalities involving $B$ yield
\[\varphi_m\big(((A_i\cap B)\cup(A_i'\sm B))\cap U\big)\le\]
\[\le\varphi_m((A_i\cap B)\cap U)+\varphi_m((A_i'\sm B)\cap U)\le3\cdot2^{-q}\lambda(U)/m.\]
Thus, for every $m>m_q$ we have (recall again that $2^{-q}<\eps/4$)
\[\varphi_m\big(((A_i\cap B)\cup(A_i'\sm B))\cap U\big)\le3\cdot2^{-q}\lambda(U)/m<\eps\lambda(U)/m.\]
It follows that the algebra $\tilde{\bB}$ satisfies condition (T.2) of property (T).

\medskip

Ultimately, we get that the algebra $\tilde{\bB}$ has property (T) and thus the proof of the lemma is finished.
\end{proof}


\section{Proof of Theorem \ref{theorem:main}\label{sec:proof}}


Theorem \ref{theorem:main} immediately follows from the following result.

\begin{theorem}\label{theorem:main_prop_t}
Assuming that $\frakp=\frakc$, there exists a Boolean subalgebra $\aA$ of $Bor(\Cantor)$ with the Grothendieck property and property (T), and hence without the Nikodym property.
\end{theorem}
\begin{proof}
We assume that $\frakp=\frakc$ or, equivalently by Bell's theorem (see \cite{Bel81}), that $\mathsf{MA}(\sigma\text{-centered})$ holds. Let us fix a bijection $f\colon\frakc\to\frakc\times\frakc\times\frakc$ such that for every $\alpha<\frakc$ if $f(\alpha)=(\alpha_1,\alpha_2,\alpha_3)$, then $\alpha\ge\alpha_1$.

\medskip

We will construct the algebra $\aA$ by induction on $\alpha<\frakc$, that is, $\aA$ will be the union of an inductively obtained increasing sequence $\seq{\bB_\alpha}{\alpha<\frakc}$ of Boolean subalgebras of $Bor(\Cantor)$. Those subalgebras will have the following properties: $\bB_0=Clopen(\Cantor)$, $\bB_\beta$ has property (T) and $|\bB_\beta|<\frakp$ for every $\beta<\frakc$,  $\bB_\beta=\bigcup_{\alpha<\beta}\bB_\alpha$ for every limit $\beta<\frakc$, and if $\beta=\alpha+1<\frakc$ and $f(\alpha)=(\alpha_1,\alpha_2,\alpha_3)$, then:
\begin{enumerate}[(i)]
	\item we fix a sequence $\seq{\cC(\alpha_1,\xi)}{\xi<\frakc}$ of countable Boolean subalgebras of $\bB_{\alpha_1}$ such that for every countable subset $C$ of $\bB_{\alpha_1}$ there is $\xi<\frakc$ with $C\subseteq\cC(\alpha_1,\xi)$ (such a sequence exists since $|\bB_{\alpha_1}|<\frakc$); we also fix: 
	\begin{itemize}
		\item an enumeration $\seq{\ol{\nu}(\alpha_1,\alpha_2,\zeta)}{\zeta<\frakc}$ of sequences of measures on $\cC(\alpha_1,\alpha_2)$,
		\item an enumeration $\seq{\ol{H}(\alpha_1,\alpha_2,\zeta)}{\zeta<\frakc}$ of countable antichains in $\cC(\alpha_1,\alpha_2)$,
		\item and an enumeration $\seq{\theta(\alpha_1,\alpha_2,\zeta)}{\zeta<\frakc}$ of finite non-negative measures on $\cC(\alpha_1,\alpha_2)$,
	\end{itemize}
	such that for every $\zeta<\frakc$ we have:
	\begin{enumerate}[(a)]
		\item $\big\|\ol{\nu}(\alpha_1,\alpha_2,\zeta)(n)\big\|=1$ for every $n\io$,
		\item $\big|\ol{\nu}(\alpha_1,\alpha_2,\zeta)(n)\big|\big(\ol{H}(\alpha_1,\alpha_2,\zeta)(n)\big)\ge0.95$ for every $n\io$,
		\item for every $A\in\cC(\alpha_1,\alpha_2)$ it holds
		\[\lim_{n\to\infty}\big|\ol{\nu}(\alpha_1,\alpha_2,\zeta)(n)\big|(A)=\theta(\alpha_1,\alpha_2,\zeta)(A),\]
	\end{enumerate}
	\item there are $B\in\bB_\beta$ and a strictly increasing sequence $\seqk{n_k\io}$ such that for every $k\io$ we have:
	\begin{enumerate}[(a)]
		\item $B\cap\ol{H}(\alpha_1,\alpha_2,\alpha_3)(n_k)\in\cC(\alpha_1,\alpha_2)$,
		\item $\big|\ol{\nu}(\alpha_1,\alpha_2,\alpha_3)(n_{2k})\big(B\cap\ol{H}(\alpha_1,\alpha_2,\alpha_3)(n_{2k})\big)\big|\le0.1$,
		\item $\big|\ol{\nu}(\alpha_1,\alpha_2,\alpha_3)(n_{2k+1})\big(B\cap\ol{H}(\alpha_1,\alpha_2,\alpha_3)(n_{2k+1})\big)\big|\ge0.3$.
	\end{enumerate}
\end{enumerate}

%

\medskip

As said, start with $\bB_0=Clopen(\Cantor)$ (so $\bB_0$ has property (T)), and assume that for some $0<\beta<\frakc$ the initial sequence $\seq{\bB_\xi}{\xi<\beta}$ of Boolean subalgebras of $Bor(\Cantor)$ with property (T), cardinality $<\frakc$, and satisfying at the successor steps conditions (i) and (ii), has been constructed. If $\beta$ is a limit ordinal, then set $\bB_\beta=\bigcup_{\xi<\beta}\bB_\xi$, so $\bB_\beta$ has property (T) by Lemma \ref{lemma:union_prop_t}. If $\beta=\alpha+1$ for some $\alpha<\frakc$, then write $f(\alpha)=(\alpha_1,\alpha_2,\alpha_3)$. For
\begin{itemize}
	\item $\bB=\bB_\alpha$, 
	\item $\cC=\cC(\alpha_1,\alpha_2)$, 
	\item $\seqn{\nu_n}$ and $\seqn{H_n}$ such that $\nu_n=\ol{\nu}(\alpha_1,\alpha_2,\alpha_3)(n)$ and $H_n=\ol{H}(\alpha_1,\alpha_2,\alpha_3)(n)$ for every $n\io$,
	\item and $\theta=\theta(\alpha_1,\alpha_2,\alpha_3)$, 
\end{itemize}
by Lemma \ref{lemma:premiere}, there exist a set $B\in Bor(\Cantor)$ and a strictly increasing sequence $\seqk{n_k\io}$ satisfying condition (ii), and such that the algebra $\tilde{\bB}$ generated by $\bB\cup\{B\}$ has property (T). Since $|\bB|<\frakp$, we also have $|\tilde{\bB}|<\frakp$. Setting $\bB_\beta=\tilde{\bB}$ finishes the inductive step.

\medskip

Proceed as described in the previous paragraph for every $\alpha<\frakc$. After the construction of the sequence $\seq{\bB_\alpha}{\alpha<\frakc}$ has been finished, put
\[\aA=\bigcup_{\alpha<\frakc}\bB_\alpha.\]
By Lemmas \ref{lemma:union_prop_t} and \ref{lemma:no_nik}, $\aA$ does not have the Nikodym property, 
so we only need to show that $\aA$ has the Grothendieck property. We will use Lemma \ref{lemma:kpr}. Let $\seqn{\nu_n}$ be a sequence of measures on $\aA$ such that $\|\nu_n\|=1$ for every $n\io$ and there exists an antichain $\seqn{H_n}$ in $\aA$ such that $|\nu_n|(H_n)\ge0.95$ for every $n\io$. Since $\cf(\frakc)\ge\omega_1$, there are $\alpha_1,\alpha_2<\frakc$ such that $\seqn{H_n}\sub\cC(\alpha_1,\alpha_2)$ and $\|\nu_n\rstr\cC(\alpha_1,\alpha_2)\|=1$ for every $n\io$ (note that to attain the norm of a measure, we only need countably many elements of a Boolean algebra). 
Since $\cC(\alpha_1,\alpha_2)$ is countable, by going to a subsequence if necessary, we can find a finite non-negative measure $\theta$ on $\cC(\alpha_1,\alpha_2)$ such that $\lim_{n\to\infty}|\nu_n|(A)=\theta(A)$ for every $A\in\cC(\alpha_1,\alpha_2)$.

There is $\alpha_3<\frakc$ such that $\theta=\theta(\alpha_1,\alpha_2,\alpha_3)$ and for every $n\io$ we have
\[\nu_n\rstr\cC(\alpha_1,\alpha_2)=\ol{\nu}(\alpha_1,\alpha_2,\alpha_3)(n)\quad\text{and}\quad H_n=\ol{H}(\alpha_1,\alpha_2,\alpha_3)(n).\]
Let $\alpha<\frakc$ be a unique ordinal number such that $f(\alpha)=(\alpha_1,\alpha_2,\alpha_3)$. Let $B\in\bB_{\alpha+1}$ and $\seqk{n_k}$ be as in property (ii). For every $k\io$ we have
\[|\nu_{n_{2k}}(B)|=\big|\nu_{n_{2k}}(B\cap H_{n_{2k}})+\nu_{n_{2k}}(B\sm H_{n_{2k}})\big|\le|\nu_{n_{2k}}(B\cap H_{n_{2k}})|+|\nu_{n_{2k}}(B\sm H_{n_{2k}})|\le\]
\[\le0.1+|\nu_{n_{2k}}|(H_{n_{2k}}^c)\le0.1+0.05=0.15,\]
as well as
\[|\nu_{n_{2k+1}}(B)|=\big|\nu_{n_{2k+1}}(B\cap H_{n_{2k+1}})+\nu_{n_{2k+1}}(B\sm H_{n_{2k+1}})\big|\ge\]
\[\ge|\nu_{n_{2k+1}}(B\cap H_{n_{2k+1}})|-|\nu_{n_{2k+1}}(B\sm H_{n_{2k+1}})|\ge\]
\[\ge0.3-|\nu_{n_{2k+1}}|(H_{n_{2k+1}}^c)\ge0.3-0.05=0.25,\]
which shows that the limit $\lim_{n\to\infty}\nu_n(B)$ does not exist. By Lemma \ref{lemma:kpr}, $\aA$ has the Grothendieck property. The proof of the theorem is thus finished.
\end{proof}

\begin{theorem}\label{theorem:main_no_prop_t}
Assuming that $\frakp=\frakc$, there exists a Boolean subalgebra $\bB$ of $Bor(\Cantor)$ with the Grothendieck property but without the Nikodym property and without property (T).
\end{theorem}
\begin{proof}
Let $A=V_0^0$. Since the Boolean algebras $Bor(\Cantor)$ and $Bor(A)$ are isomorphic, Theorem \ref{theorem:main_prop_t} implies that there is a Boolean subalgebra $\aA$ of $Bor(A)$ with the Grothendieck property but without the Nikodym property. Define a Boolean subalgebra $\bB$ of $Bor(\Cantor)$ in the following way:
\[\bB=\big\{B\cup C\colon\  B\in\aA, C\in Bor(\Cantor\sm A)\big\}.\]
Since $Bor(\Cantor\sm A)\sub\bB$, by Proposition \ref{prop:open_sets_prop_t} (applied to $W=\Cantor\sm A=V_0^1$), $\bB$ does not have property (T). Moreover, since both $\aA$ and $Bor(\Cantor\sm A)$ have the Grothendieck property (note that $Bor(\Cantor\sm A)$ is $\sigma$-complete), $\bB$ has the Grothendieck property, too. On the other hand, as $\aA$ does not have the Nikodym property, $\bB$ does not have the Nikodym property either.
\end{proof}

\section{Open questions\label{sec:questions}}

In this section we present several open questions concerning Boolean algebras with the Grothendieck property but without the Nikodym property. We start with the most important one.

\begin{question}\label{ques:zfc}
Does there exist in ZFC a Boolean algebra with the Grothendieck property but without the Nikodym property?
\end{question}

Based on constructions presented in this paper and in \cite{Tal80}, it seems that a careful analysis of property (T) and its general relations to Boolean algebras with the Grothendieck property but without the Nikodym property could pave the way for making an essential step towards answering Question \ref{ques:zfc}. We therefore pose below several auxiliary questions concerning property (T).

The Boolean algebras constructed in Theorem \ref{theorem:main_prop_t} as well as in \cite{Tal80} both have property (T). The Boolean algebra $\bB$ from Theorem \ref{theorem:main_no_prop_t} does not have property (T), but the reason why $\bB$ does not have the Nikodym property either is that there exists an element $A\in\bB$ such that the Boolean subalgebra $\{B\cup(\Cantor\sm A)\in\bB\colon B\le A\}$ of $\bB$  has property (T). It seems thus interesting whether property (T) (or its variant of some kind) is indeed in any sense necessary for a Boolean subalgebra of $Bor(\Cantor)$ to have the Grothendieck property but not the Nikodym property. We thus ask the following question.


\begin{question}\label{ques:t_con_gr}
Is it consistent that every Boolean subalgebra of $Bor(\Cantor)$ with property (T) does not have the Grothendieck property?
\end{question}

The next few problems concern maximal (in the sense of inclusion) Boolean subalgebras of $Bor(\Cantor)$ with property (T).

\begin{question}\label{ques:t_max_gr}
Assume that $\aA$ is a maximal Boolean subalgebra of $Bor(\Cantor)$ with property (T). Does $\aA$ necessarily have the Grothendieck property?
\end{question}

Since $Clopen(\Cantor)$ has property (T) and the union of an increasing chain of Boolean subalgebras of $Bor(\Cantor)$ with property (T) has property (T), an affirmative answer (in ZFC) to Question \ref{ques:t_max_gr} would also yield a positive answer to Question \ref{ques:zfc}.

Let us define the following cardinal characteristic of the continuum:
\[\mathfrak{prop}_{(T)}=\min\big\{|\aA|\colon\ \aA\text{ is a maximal subalgebra of }Bor(\Cantor)\text{ with property (T)}\big\}.\]
Using methods from Section \ref{sec:ext_prop_t} one can show that $\omega_1\le\mathfrak{prop}_{(T)}\le\frakc$. The Continuum Hypothesis or Martin's axiom (or its variant for $\sigma$-centered posets) imply that $\mathfrak{prop}_{(T)}=\frakc$. We do not now whether consistently it may hold that $\mathfrak{prop}_{(T)}<\frakc$.

\begin{question}\label{ques:propt_less_than_c}
Is there a model of ZFC in which the inequality $\mathfrak{prop}_{(T)}<\frakc$ holds?
\end{question}

If such a model exists, then one can naturally ask whether $\mathfrak{prop}_{(T)}$ is already equal to one of many well-known cardinal characteristics of the continuum, e.g., to one occurring in classical van Douwen's diagram or Cicho\'{n}'s diagram (see \cite{Bla10} or \cite{vD84} for information on these standard diagrams). Note that the proof of Theorem \ref{theorem:main} actually shows that $\mathfrak{prop}_{(T)}\ge\frakp$. 

\begin{question}\label{ques:propt_standard}
Is $\mathfrak{prop}_{(T)}$ equal in ZFC to some standard cardinal characteristic of the continuum, e.g., to one occurring in van Douwen's diagram or Cicho\'{n}'s diagram? 

If not, what are lower and upper bounds for $\mathfrak{prop}_{(T)}$ in terms of standard cardinal characteristics?
\end{question}

If Question \ref{ques:t_max_gr} has an affirmative answer, then $\mathfrak{prop}_{(T)}\ge\frakgr$, where $\frakgr$ denotes the minimal cardinality of an infinite Boolean algebra with the Grothendieck property. The number $\frakgr$ was studied in \cite{SZforext} and \cite{SZadding}, see also the survey \cite{SobKok}. One can show e.g. that $\frakgr\ge\max(\cov(\mM),\fraks)$, where $\cov(\mM)$ denotes the covering number of the meager $\sigma$-ideal $\mM$ and $\fraks$ denotes the splitting number (see \cite{Bla10} for basic information on $\cov(\mM)$ and $\fraks$).

\begin{question}
Is it consistent that $\mathfrak{prop}_{(T)}<\frakgr$?
\end{question}

\section*{Appendix. Proof of Lemma \ref{lemma:combinatorial}}

The following proof of Lemma \ref{lemma:combinatorial} is due to Talagrand \cite[p. 170--171]{Tal84}. Since Lemma \ref{lemma:combinatorial} is crucial for our paper as well as Talagrand's original presentation of the argument might be considered a bit too concise, we present the proof of the lemma in a detailed way.

We will need the following variant of the Paley--Zygmund theorem, being essentially \cite[Theorem 3, page 31]{Kah94}. For a random variable $X$ on a probability space $(\Omega,\fF,\pP)$ we denote by $\eE(X)$ and $V(X)$ the expected value and the variance of $X$, respectively.

\begin{theorem}[Paley--Zygmund]\label{theorem:paleyzygmund}
Fix $l\io$ and $C>0$. Suppose that $\{X_i\colon 1\le i\le l\}$ is a family of real-valued random variables on a probability space $(\Omega,\fF,\pP)$ such that $X_i\in L_4(\Omega,\fF,\pP)$, $\eE(X_i)=0$, and $\eE(X_i^4)\le CV(X_i)^2$ for all $1\le i\le l$. Let also $\xi\in(0,1)$. Then,
\[\pP\big(\big\{y\in\Omega\colon\ |X_1(y)+\ldots+X_l(y)|>\xi(V(X_1)+\ldots+V(X_l))^{1/2}\big\}\big)>\zeta,\]
where $\zeta=(1-\xi^2)^2\cdot\min\{1/3,1/C\}$.
\end{theorem}

\begin{proof}[{\textbf{Proof of Lemma \ref{lemma:combinatorial}}}]
Let $n_0\io$, $n_0\ge t$, be such that $n^3\cdot2^{-n}\le\eta$ for every $n\ge n_0$. It follows that
\[\tag{L.1}\big(\eta n2^{-n}\big)^{1/2}\le\eta/n\]
for every $n\ge n_0$. Fix $n\ge n_0$ and $Q\in\aA_n$ such that there exists an atom $G$ of $\aA_t$ for which we have $\lambda(Q\cap G)\ge0.95\lambda(G)$. Without loss of generality let us assume that $Z\cup R\sub Q$, so $Z\cap Q=Z$ and $R\cap Q=R$. Set $Q'=Q\sm(Z\cup R\cup P)$. Of course, $Q'\in\aA_n$. We also have $\lambda(Q'\cap G)\ge0.92\lambda(G)$. Indeed, note that
\[\lambda(Z\cup R\cup P)\le\lambda(Z)+\lambda(R)+\lambda(P)\le3\cdot2^{-8}\eta^2<3\cdot2^{-8}\cdot2^{-10}\cdot2^{-t-1}<2^{-16}\lambda(G)<0.03\lambda(G),\]
and so
\[\lambda(Q'\cap G)=\lambda(Q\cap G)-\lambda(Q\cap(Z\cup R\cup P)\cap G)\ge\lambda(Q\cap G)-\lambda(Z\cup R\cup P)>\]
\[>0.95\lambda(G)-0.03\lambda(G)=0.92\lambda(G).\]
Consequently, we have
\[\lambda(Q')\ge\lambda(Q'\cap G)\ge0.92\lambda(G)=0.92\cdot2^{-t-1}>2^{-t-11}>\eta.\]

For every $t<m\le n$ set
\[a_m=\lambda(Z\cap V_m^0)-\lambda(Z\cap V_m^1)-\lambda(R\cap V_m^0)+\lambda(R\cap V_m^1),\]
so $|a_m|=\psi_m(Z,R)$, and, additionally, for every subset $N\sub Q'$ define
\[b_m(N)=\lambda(N\cap V_m^0)-\lambda(N\cap V_m^1),\]
so $|b_m(N)|=\varphi_m(N)$. Let $\eta'\in(\eta/2,\eta)$ be of the form $\eta'=k2^{-(n+1)}$ for some $k\io$ (recall that $n^3\cdot2^{-n}\le\eta<2^{-t-11}$). 
Since $\lambda(Q')>\eta$, there is $M\sub Q'$ such that $M\in\aA_n$, $\lambda(M)=\eta'$, and the number
\[S(M)=\sum_{t<m\le n}(b_m(M)+a_m)^2\]
is minimal. For every $t<m\le n$ write shortly $b_m=b_m(M)$. Since $M\cap Z=\emptyset$, we have
\[|a_m+b_m|=\]
\[=\big|\lambda(Z\cap V_m^0)-\lambda(Z\cap V_m^1)-\lambda(R\cap V_m^0)+\lambda(R\cap V_m^1)+\lambda(M\cap V_m^0)-\lambda(M\cap V_m^1)\big|=\]
\[=\big|\lambda((M\cup Z)\cap V_m^0)-\lambda((M\cup Z)\cap V_m^1)-\lambda(R\cap V_m^0)+\lambda(R\cap V_m^1)\big|=\]
\[=\psi_m(M\cup Z, R)=\psi_m(M\cup(Z\cap Q),\ R\cap Q),\]
so to finish the proof of the lemma we need to show that
\[\tag{L.2}|a_m+b_m|\le\eta/m.\]

\medskip

For $t<m\le n$ and an atom $W$ of $\aA_n$ define
\[x_m^W=\begin{cases}
1,&\text{ if }W\sub V_m^0,\\
-1,&\text{ if }W\sub V_m^1.
\end{cases}\]
Since $\lambda(M)<\eta<\lambda(Q')$, we have $Q'\sm M\neq\emptyset$. Also, $Q'\sm M\in\aA_n$. Let $U$ and $V$ be any atoms of $\aA_n$ such that $U\sub M$ and $V\sub Q'\sm M$. Put
\[M^{U,V}=(M\sm U)\cup V.\]
Since $U\sub M$ and $M\cap V=\emptyset$, for $t<m\le n$ we have
\[b_m(M^{U,V})=\lambda(M^{U,V}\cap V_m^0)-\lambda(M^{U,V}\cap V_m^1)=\]
\[=\lambda\big(((M\sm U)\cup V)\cap V_m^0\big)-\lambda\big(((M\sm U)\cup V)\cap V_m^1\big)=\]
\[=\lambda(M\cap V_m^0)-\lambda(M\cap V_m^1)-\lambda(U\cap V_m^0)+\lambda(U\cap V_m^1)+\lambda(V\cap V_m^0)-\lambda(V\cap V_m^1)=\]
\[=b_m+2^{-n-1}(x_m^V-x_m^U).\]
Note that $\lambda(M^{U,V})=\eta'$. By the choice of $M$, we have
\[0\le S(M^{U,V})-S(M)=\sum_{t<m\le n}(b_m(M^{U,V})+a_m)^2-\sum_{t<m\le n}(b_m+a_m)^2=\]
\[=\sum_{t<m\le n}\big(b_m+2^{-n-1}(x_m^V-x_m^U)+a_m\big)^2-\sum_{t<m\le n}(b_m+a_m)^2=\]
\[=\sum_{t<m\le n}2(b_m+a_m)\cdot 2^{-n-1}(x_m^V-x_m^U)+\sum_{t<m\le n}2^{-2n-2}(x_m^V-x_m^U)^2.\]
Multiplying the above inequality by $2^{n}$, we get
\[\sum_{t<m\le n}(b_m+a_m)(x_m^V-x_m^U)+2^{-n-2}\sum_{t<m\le n}(x_m^V-x_m^U)^2\ge0.\]
Note that $(x_m^V-x_m^U)^2\le4$ for every $t<m\le n$, so it holds
\[\tag{L.3}\sum_{t<m\le n}(b_m+a_m)(x_m^V-x_m^U)+n2^{-n}\ge0.\]

Let the mapping $T\colon\Cantor\to\Cantor$ be defined for every $z\in\Cantor$ and $m\io$ as follows:
\[T(z)(m)=\begin{cases}
z(m),&\text{ if }m\not\in(t,n],\\
1-z(m),&\text{ if }m\in(t,n].
\end{cases}\]
Then, $T$ is a measure-preserving homeomorphism of the Cantor space $\Cantor$. Since
\[\lambda(M)<\eta<2^{-t-11}=2^{-10}\lambda(G)<0.001\lambda(G),\]
we have
\[\lambda((Q'\sm M)\cap G)=\lambda(Q'\cap G)-\lambda(Q'\cap M\cap G)\ge\]
\[\ge\lambda(Q'\cap G)-\lambda(M)>0.92\lambda(G)-0.001\lambda(G)>0.9\lambda(G).\]
As $G$ is an atom of $\aA_t$, it also holds
\[\lambda(T[Q'\sm M]\cap G)=\lambda(T[(Q'\sm M)\cap G])=\lambda((Q'\sm M)\cap G)>0.9\lambda(G),\]
so
\[\lambda(T[Q'\sm M]^c\cap G)<0.1\lambda(G).\]
Consequently, we get
\[\lambda\big((Q'\sm M)\cap T[Q'\sm M]\big)\ge\lambda\big((Q'\sm M)\cap T[Q'\sm M]\cap G\big)=\]
\[=\lambda((Q'\sm M)\cap G)-\lambda\big((Q'\sm M)\cap T[Q'\sm M]^c\cap G\big)\ge\]
\[\ge\lambda((Q'\sm M)\cap G)-\lambda(T[Q'\sm M]^c\cap G)\ge\]
\[\ge0.9\lambda(G)-0.1\lambda(G)=0.8\lambda(G),\]
so in particular $(Q'\sm M)\cap T[Q'\sm M]$ is a non-empty set of positive $\lambda$-measure.

We will now use Theorem \ref{theorem:paleyzygmund}. For each $m\io$ let the function $\pi_m\colon\Cantor\to\{1,-1\}$ be defined for every $z\in\Cantor$ as follows: $\pi_m(z)=1$ if and only if $z(m)=0$ (so, in short, $\pi_m(z)=(-1)^{z(m)}$). Let $\Omega=(Q'\sm M)\cap T[Q'\sm M]$ (endowed with the subspace topology of $\Cantor$), $\fF=Bor(\Omega)$, $\pP=\lambda/\lambda(\Omega)$, $\xi=2^{-2}$, $l=n-t$, and for every $1\le i\le l$ let the random variable $X_i\colon\Omega\to\R$ be such that $X_i(y)=\pi_{t+i}(y)(b_{t+i}+a_{t+i})$ for all $y\in\Omega$ ($\sub\Cantor$). Fix $1\le i\le l$ for a moment. Since $\Omega=T[\Omega]$, it is easy to check that $\eE(X_i)=0$. We have $X_i^2(y)=(b_{t+i}+a_{t+i})^2$ for all $y\in\Omega$, so $V(X_i)=(b_{t+i}+a_{t+i})^2$. Similarly, $\eE(X_i^4)=(b_{t+i}+a_{t+i})^4$. Thus, $\eE(X_i^4)\le CV(X_i)^2$ for $C=1$.

Applying Theorem \ref{theorem:paleyzygmund} for the above setting, we conclude that
\[\pP\big(\Big\{y\in\Omega\colon\ \big|\sum_{i=1}^{n-t}\pi_{t+i}(y)(b_{t+i}+a_{t+i})\big|>2^{-2}S(M)^{1/2}\Big\}\big)=\]
\[\pP\big(\Big\{y\in\Omega\colon\ \big|\sum_{i=1}^{n-t}\pi_{t+i}(y)(b_{t+i}+a_{t+i})\big|>2^{-2}\Big(\sum_{i=1}^{n-t}(b_{t+i}+a_{t+i})^2\Big)^{1/2}\Big\}\big)=\]
\[\pP\big(\big\{y\in\Omega\colon\ |X_1(y)+\ldots+X_l(y)|>\xi(V(X_1)+\ldots+V(X_l))^{1/2}\big\}\big)>\]
\[>(1/3)\cdot(1-1/4^2)^2>0.\]
It follows that there exists $y\in\Omega$ such that
\[\big|\sum_{i=1}^{n-t}\pi_{t+i}(y)(b_{t+i}+a_{t+i})\big|>2^{-2}S(M)^{1/2}.\]
Let $V$ be an atom of $\aA_n$ such that $y\in V$. Note that $V\sub\Omega$ (as $\Omega\in\aA_n$). Since $\pi_{t+i}(y)=x_{t+i}^V$ for every $1\le i\le n-t$, we get
%
\[\big|\sum_{i=1}^{n-t}x_{t+i}^V(b_{t+i}+a_{t+i})\big|=\big|\sum_{t<m\le n}x_m^V(b_m+a_m)\big|\ge 2^{-2}S(M)^{1/2}.\]
Replacing $y$ with $T(y)$ and $V$ with $T[V]$, if necessary, we can assume that
\[\sum_{t<m\le n}x_m^V(b_m+a_m)\le -2^{-2}S(M)^{1/2}.\]
Consequently, by inequality (L.3), for any atom $U$ of $\aA_n$ such that $U\sub M$ we have
\[\tag{L.4}\sum_{t<m\le n}x_m^U(b_m+a_m)\le\sum_{t<m\le n}x_m^V(b_m+a_m)+n2^{-n}\le n2^{-n}-2^{-2}S(M)^{1/2}.\]

Since for every $t<m\le n$ it holds
\[x_m^U=2^{n+1}(\lambda(U\cap V_m^0)-\lambda(U\cap V_m^1)),\]
we have
\[2^{n+1}\sum_{t<m\le n}b_m(b_m+a_m)=2^{n+1}\sum_{t<m\le n}\big(\lambda(M\cap V_m^0)-\lambda(M\cap V_m^1)\big)\cdot(b_m+a_m)=\]
\[=2^{n+1}\sum_{t<m\le n}\sum_{\substack{U\sub M\\U\text{ -- atom of }\aA_n}}\big(\lambda(U\cap V_m^0)-\lambda(U\cap V_m^1)\big)\cdot(b_m+a_m)=\]
\[=\sum_{\substack{U\sub M\\U\text{ -- atom of }\aA_n}}\sum_{t<m\le n}2^{n+1}\big(\lambda(U\cap V_m^0)-\lambda(U\cap V_m^1)\big)\cdot(b_m+a_m)=\]
\[=\sum_{\substack{U\sub M\\U\text{ -- atom of }\aA_n}}\sum_{t<m\le n}x_m^U(b_m+a_m)\le\]
\[\le\big|\{U\colon\ U\text{ -- atom of }\aA_n,\ U\sub M\}\big|\cdot\big(n2^{-n}-2^{-2}S(M)^{1/2}\big)=\]
\[=k\big(n2^{-n}-2^{-2}S(M)^{1/2}\big),\]
where the last inequality follows from (L.4) and the last equality is a consequence of the fact that $\lambda(M)=\eta'=k2^{-(n+1)}$. Dividing the above inequality by $2^{n+1}$, we get
\[\tag{L.5}\sum_{t<m\le n}b_m(b_m+a_m)\le k2^{-(n+1)}\big(n2^{-n}-2^{-2}S(M)^{1/2}\big)=\eta'\big(n2^{-n}-2^{-2}S(M)^{1/2}\big).\]

By the Cauchy--Schwarz inequality, we have
\[\Big(\Big(\sum_{t<m\le n}a_m^2\Big)^{1/2}\cdot S(M)^{1/2}\Big)^2=\Big(\sum_{t<m\le n}a_m^2\Big)\cdot S(M)=\Big(\sum_{t<m\le n}a_m^2\Big)\cdot\Big(\sum_{t<m\le n}(b_m+a_m)^2\Big)\ge\]
\[\ge\Big(\sum_{t<m\le n}a_m(b_m+a_m)\Big)^2.\]
Taking the square root, we obtain
\[\tag{L.6}\Big(\sum_{t<m\le n}a_m^2\Big)^{1/2}\cdot S(M)^{1/2}\ge\sum_{t<m\le n}a_m(b_m+a_m).\]

By inequalities (L.5) and (L.6), we have
\[S(M)=\sum_{t<m\le n}(b_m+a_m)^2=\]
\[=\sum_{t<m\le n}b_m(b_m+a_m)+\sum_{t<m\le n}a_m(b_m+a_m)\le\]
\[\le\eta'\big(n2^{-n}-2^{-2}S(M)^{1/2}\big)+\Big(\sum_{t<m\le n}a_m^2\Big)^{1/2}\cdot S(M)^{1/2}=\]
\[=\eta'n2^{-n}+S(M)^{1/2}\big(-\eta'/4+\Big(\sum_{t<m\le n}a_m^2\Big)^{1/2}\big),\]
or, in short,
\[\tag{L.7}S(M)\le\eta'n2^{-n}+S(M)^{1/2}\big(-\eta'/4+\Big(\sum_{t<m\le n}a_m^2\Big)^{1/2}\big).\]

We will now estimate the value of $\big(\sum_{t<m\le n}a_m^2\big)^{1/2}$. Note that for any set $N\in Bor(\Cantor)$ we have
\[\lambda(Z\cap N)-\lambda(R\cap N)=\lambda((Z\sm R)\cap N)-\lambda((R\sm Z)\cap N),\]
so
\[\Big(\sum_{t<m\le n}a_m^2\Big)^{1/2}=\Big(\sum_{t<m\le n}\big(\lambda(Z\cap V_m^0)-\lambda(Z\cap V_m^1)-\lambda(R\cap V_m^0)+\lambda(R\cap V_m^1)\big)^2\Big)^{1/2}=\]
\[=\Big(\sum_{t<m\le n}\big(\lambda((Z\sm R)\cap V_m^0)-\lambda((Z\sm R)\cap V_m^1)-\lambda((R\sm Z)\cap V_m^0)+\lambda((R\sm Z)\cap V_m^1)\big)^2\Big)^{1/2}=\]
\[=\Big(\sum_{t<m\le n}\Big(\int_{\Cantor}\chi_{Z\sm R}\pi_md\lambda-\int_{\Cantor}\chi_{R\sm Z}\pi_md\lambda\Big)^2\Big)^{1/2}.\]
Note that the sequence $\seqm{\pi_m}$ constitutes an orthonormal system in the real Hilbert space $L_2\big(\Cantor,Bor(\Cantor),\lambda\big)$, endowed with the inner product $\langle f,g\rangle=\int_{\Cantor}fgd\lambda$. By Bessel's inequality, we have
\[\Big(\sum_{t<m\le n}\Big(\int_{\Cantor}\chi_{Z\sm R}\pi_md\lambda-\int_{\Cantor}\chi_{R\sm Z}\pi_md\lambda\Big)^2\Big)^{1/2}=\Big(\sum_{t<m\le n}\Big(\int_{\Cantor}(\chi_{Z\sm R}-\chi_{R\sm Z})\pi_md\lambda\Big)^2\Big)^{1/2}=\]
\[=\Big(\sum_{t<m\le n}\big\langle\chi_{Z\sm R}-\chi_{R\sm Z},\ \pi_m\big\rangle^2\Big)^{1/2}=\Big(\sum_{t<m\le n}\big|\big\langle\chi_{Z\sm R}-\chi_{R\sm Z},\ \pi_m\big\rangle\big|^2\Big)^{1/2}\le\]
\[\le\|\chi_{Z\sm R}-\chi_{R\sm Z}\|_{L_2(\Cantor)}=\Big(\int_{\Cantor}(\chi_{Z\sm R}-\chi_{R\sm Z})^2d\lambda\Big)^{1/2}\le\Big(\int_{\Cantor}\chi_{Z\cup R}^2d\lambda\Big)^{1/2}=\lambda(Z\cup R)^{1/2}.\]
Since $\eta/2<\eta'<\eta$, we thus obtain
\[\Big(\sum_{t<m\le n}a_m^2\Big)^{1/2}\le\lambda(Z\cup R)^{1/2}\le(\lambda(Z)+\lambda(R))^{1/2}\le(2\cdot2^{-8}\cdot\eta^2)^{1/2}=2^{-7/2}\eta<\]
\[<2^{-3}\eta=\eta/8<\eta'/4,\]
so
\[\tag{L.8}\Big(\sum_{t<m\le n}a_m^2\Big)^{1/2}<\eta'/4.\]
Applying (L.8) to (L.7), we get
\[S(M)\le\eta'n2^{-n}\le\eta n2^{-n},\]
which, by (L.1), yields for every $t<m\le n$ that
\[|b_m+a_m|\le\Big(\sum_{t<m'\le n}(b_{m'}+a_{m'})^2\Big)^{1/2}=S(M)^{1/2}\le\eta/n\le\eta/m.\]
This proves inequality (L.2) and the proof of the lemma is thus finished.
\end{proof}

\end{document}